\documentclass[12pt]{amsart}
\usepackage{amsmath}
\usepackage{amssymb}

\numberwithin{equation}{section}

\newtheorem{defn}{Definition}[section]
\newtheorem{theorem}{Theorem}[section]

\newtheorem{conjecture}[theorem]{Conjecture}
\newtheorem{corollary}[theorem]{Corollary}

\newtheorem{lemma}[theorem]{Lemma}
\newtheorem{prop}[theorem]{Proposition}
\newtheorem{remark}[theorem]{Remark}
\def\begineq{\begin{equation}}
\def\endeq{\end{equation}}

\def \n{\noindent}

\def \v{\vskip 0.1in}

\def \mc{\mathcal}

\def \mk{\mathfrak}

\def \cplane{\mathbb{C}}
\def \rone{\mathbb{R}}
\def \pone{\mathbb{P}}
\def \integer{\mathbb{Z}}

\def \mfp{\mathfrak{p}^s}
\def \mfq{\mathfrak{q}^s}
\def \mfP{\mathfrak{p}^{sf}}
\def \mfQ{\mathfrak{q}^{sf}}

\def \om{\overline{\mathcal{M}}}
\def \M{\mathcal{M}}

\title{Singular symplectic flops and Ruan cohomology}
\author{Bohui Chen}
\address{Department of Mathematics, Sichuan University,
        Chengdu,610064, China}
\email{bohui@cs.wisc.edu}
\author{An-Min Li}
\address{Department of Mathematics, Sichuan University,
        Chengdu,610064, China}
\email{math$\_$li@yahoo.com.cn}
\author{Qi Zhang}
\address{Department of Mathematics, University of Missouri-Columbia}
\email{qi@math.missouri.edu}
\author{Guosong Zhao }
\address{Department of Mathematics, Sichuan University,
        Chengdu,610064, China}
\email{gszhao@scu.edu.cn}
\thanks{B.C. and A.L. are supported by  NSFC,
G.Z. is supported by a grant of NSFC and Qiushi Funding.}
\date{}
\begin{document}
\maketitle \abstract In this paper, we  study the symplectic
geometry of singular conifolds of the finite group quotient
$$
W_r=\{(x,y,z,t)|xy-z^{2r}+t^2=0 \}/\mu_r(a,-a,1,0), r\geq 1,
$$
which we call  orbi-conifolds. The related orbifold symplectic
conifold transition and orbifold symplectic flops are constructed.
Let $X$ and $Y$ be two symplectic orbifolds connected by such a
flop. We study orbifold Gromov-Witten invariants of exceptional
classes on $X$ and $Y$ and show that they have isomorphic Ruan
cohomologies. Hence, we verify a conjecture of Ruan.
\endabstract
\def \rroot{e^{\frac{2\pi i}{r}}}

\section{Introduction}\label{sect_0}

In \cite{LR}, the authors proved an elegant result that any two
smooth minimal models in dimension three have the same quantum
cohomology. Besides the key role of the relative invariants
introduced in the paper, one of the main building blocks towards
this result is the understanding of  how the Gromov-Witten
invariants change under flops. The description of a smooth flop is
closely related to the conifold singularity
$$
W_1=\{(x,y,z,t)|xy-z^2+t^2=0 \}.
$$
 A crucial step in their proof is a symplectic description of a flop
 and hence symplectic techniques
 can be applied. However, it is well-known
 that the appropriate category for birational geometry is
 singular manifolds
  with terminal singularities. In  complex dimension three, terminal singularities
 are deformations of orbifolds.
  In this paper and its sequel, we initiate a program to study the
 quantum cohomology under birational transformation of orbifolds.

In the singular category,
$$
W_r=\{(x,y,z,t)|xy -z^{2r}+t^2=0\}/\mu_r(a,-a,1,0).
$$
 is a natural  replacement for the smooth conifold.
The orbifold symplectic flops coming from this model are defined
in the first part of the paper (cf. \S\ref{sect_3}). In the second part
of the paper, we compute the 3-point function of (partial) orbifold
Gromov-Witten invariants. This enables us to verify a conjecture
by Ruan in the current set-up:
 for any two symplectic orbifolds $X$ and $Y$  connected via
orbifold symplectic flops,  their Ruan cohomology rings are
isomorphic.

\subsection{Orbifold symplectic flops}\label{sect_0.1}
The singularity given by
$W_1$ has been studied intensively.
Let $\omega^o$ be the symplectic form on $W_1\setminus\{0\}$
induced from that of $\cplane^4$.
It has two small resolutions,
denoted by $W_1^s$ and $W_1^{sf}$, and a smoothing via deformation
which is denoted by $Q_1$. The transformations
$$
W^s_1\leftrightarrow Q_1,\;\;
W^{sf}_1\leftrightarrow Q_1
$$
are called {\em conifold transitions}.
And the transformation
$$
W^s_1\leftrightarrow W^{sf}_1
$$
is called a {\em  flop}.

A symplectic conifold(\cite{STY}) $(Z,\omega)$ is a space with conifold singularities
$$
P=\{p_1,\ldots\}
$$
such that $(Z\setminus P,\omega)$ is a symplectic manifold
and $\omega$ coincides with $\omega^o$ locally at $p_i\in P$. Now
suppose that $Z$ is compact and $|P|=\kappa<\infty$.
Such  $Z$ admits a smoothing, denoted by $X$, and $2^\kappa$
resolutions
$$
\mc Y=\{Y_1,\ldots, Y_{2^\kappa}\}.
$$
 In $X$  each $p_i$ is replaced
 by an exceptional sphere $L_i\cong S^3$, while for each $Y_j$, $p_i$ is replaced by
an extremal ray $\pone^1$.

In \cite{STY}, they studied a necessary and sufficient condition for
the existence of a symplectic structure on one of the $Y$ in $\mc Y$ in terms
of certain topological condition on $X$.
 They showed that {\em one of the $2^\kappa$
small resolutions
admits a symplectic stucture if and only if  on
 $X$ we have  the following  homology relation
\begin{equation}\label{eqn_1.1}
[\sum_{i=1}^\kappa
\lambda_iL_{i}]=0\in  H_3(X,\integer)
\;\;with \;\lambda_i\not= 0\; for \;all\;i.
\end{equation}
Here the $L_i$ are exceptional spheres on $X$.}

One can rephrase their theorem using cohomological language. Then,
equation \eqref{eqn_1.1} reads as
\begin{equation}\label{eqn_1.2}
[\sum_{i=1}^\kappa
\lambda_i\Theta_{i}]=0\in  H^3(X,\integer)
\;\;with \;\lambda_i\not= 0\; for \;all\;i.
\end{equation}
Here $\Theta_i$ is the Thom form of the normal
bundle of $L_i$.

The cohomological version will be generalized to the general model
with finite group quotient. Our model is
\begin{equation}\label{eqn_1.3}
W_r=\{(x,y,z,t)|xy-z^{2r}+t^2=0 \}/\mu_r(a,-a,1,0), r\geq 1.
\end{equation}
(see \cite{K} and \cite{Reid} for references). Such a local model
is called {\em $r$-conifold or an orbi-conifold} in our paper.
Such (terminal) singularities appear naturally in the Minimal Model Program.
They are the
simplest examples in the list of singularities in \cite{K}. $W_r$
without the finite quotient has been considered in \cite{Laufer}.
It also has two resolutions $\tilde W_r^{s}$ and $\tilde
W_{r}^{sf}$. We can take quotients
$$
W_r^s=\tilde W^s_r/\mu_r,\;\; W_r^{sf}=\tilde W^{sf}_r/\mu_r.
$$
 Both of them are orbifolds. In this paper, we  propose a
smoothing $Q_r$ as well. The transformations
$$
W^s_r\leftrightarrow Q_r,\;\;
W^{sf}_r\leftrightarrow Q_1
$$
are called {\em (orbi)-conifold transitions}. And the transformation
$$
W^s\leftrightarrow W^{sf}
$$
is called a {\em (orbi)-flop}.

We are interested in symplectic geometry of the orbi-conifold
$(Z,\omega_Z)$. It has a smoothing $X$ and $2^\kappa$ small
resolutions
$$
\mc Y=\{Y_i,
1\leq i\leq 2^\kappa.\}
$$
A theorem generalizing that of Smith-Thomas-Yau is
\begin{theorem}\label{theorem_1.1}
One of the $2^\kappa$
small resolutions
admits a symplectic stucture if and only if  on
 $X$ we have  the following  cohomology relation
\begin{equation}\label{eqn_1.4}
[\sum_{i=1}^\kappa
\lambda_i\Theta_{r_i}]=0\in  H^3(X,\rone)
\;\;with \;\lambda_i\not= 0\; for \;all\;i.
\end{equation}
\end{theorem}
As a corollary of this theorem, we show that if one of
$Y_i\in \mc Y$ is symplectic then so is its flop $Y^f_i\in \mc Y$ (refer to \S\ref{sect_3.1} for the definition).

\subsection{The ring structures  and Ruan's conjecture}\label{sect_0.2}
Let $X$ be an orbifold. It is well known that $H^*(X)$ does not suffice
for quantum cohomology. One should consider the so-called
twisted sectors $X_{(g)}$ on $X$ and study a bigger space
$$
H^\ast_{CR}:=H^\ast(X)\oplus\bigoplus_{(g)|g\not=1}
H^\ast(X_{(g)}).
$$
Using the orbifold Gromov-Witten invariants \cite{CR2}, one can
define the orbifold quantum ring $QH^*_{CR}(X)$.  The analogue  of
classical cohomology is known as the Chen-Ruan orbifold cohomology
ring.

Motivated by the work of Li-Ruan (\cite{LR}) on the transformation
of the quantum cohomology rings with respect to a smooth flop,
 we may ask how the
orbifold quantum cohomology ring transforms (or even how the
orbifold Gromov-Witten invariants change) via orbi-conifold
transitions or orbifold flops. It can be formulated as the
following conjecture
\begin{conjecture}\label{conjecture_1.1}
Let $Y$ be the orbifold symplectic flop of $X$, then
$$
QH^\ast_{CR}(X)\cong QH^\ast_{CR}(Y).
$$
\end{conjecture}
To completely answer the question, one  needs a full package of
technique, such as relative orbifold Gromov-Witten invariants and
degeneration formulae. These techniques are out of reach at this
moment and will be studied in  future papers(\cite{CLZZ}).

On the other hand, it is easy to show that
$$
H^\ast_{CR}(X)\cong H^\ast_{CR}(Y)
$$
additively. In general, they will have different ring structures.
In this paper, we study a new ring structure that it is in a sense
between $H^\ast_{CR}$ and $QH^\ast_{CR}$. It was first introduced
by Ruan \cite{R} in the smooth case and can be naturally extended to
orbifolds. Let's review the construction.  Let $\Gamma^s_i,
\Gamma^{sf}_i, 1\leq i\leq \kappa$ be extremal rays in $X$ and $Y$
respectively. On $X$, (and on $Y$), we use only moduli spaces of
J-curves representing multiples of $\Gamma_i$'s and define
 3-point functions on $H^\ast_{CR}(X)$ by
\begin{equation}\label{eqn_1.5}
\Psi^X_{qc}(\beta_1,\beta_2,\beta_3)=\Psi^X_{d=0}(\beta_1,\beta_2,\beta_3)
+\sum_{i=1}^\kappa\sum_{d=1}^\infty\Psi^X_{(d[\Gamma^s],0,3)}(\beta_1,\beta_2,
\beta_3).
\end{equation}
Such  functions also yield a product on
$H^\ast_{CR}(X)$. This ring is called the Ruan cohomology ring
\cite{HZ} and  denoted by $RH^\ast_{CR}(X)$.  Ruan conjectures that
{\em  if $X$, $Y$ are K-equivalent, $RH^\ast_{CR}(X)$ is
isomorphic to $RH^\ast_{CR}(Y)$.}

Our second theorem is

\begin{theorem}\label{theorem_1.2}
Suppose that  $X$ and $Y$ are connected by a sequence of
symplectic flops constructed out of $r$-conifolds. Then
$RH^\ast_{CR}(X)$ is isomorphic to $RH^\ast_{CR}(Y)$. Hence,
Ruan's conjecture holds in this case.
\end{theorem}

\v
{\em Acknowledge. }We would like to thank  Yongbin Ruan
for telling us about the program and
 for many valuable discussions. We also wish to thank
Qi Zhang, Shengda Hu and Quan Zheng for many discussions.
The second and third
authors also would like to thank  University of Wisconsin-
Madison and MSRI for their hospitality.

\section{Local Models}\label{sect_1}

\subsection{Local $r$-orbi-conifolds}\label{sect_1.1}
Let
$$
\mu_r=\langle\xi\rangle, \xi=\rroot
$$
be the cyclic group of $r$-th roots of 1. We denote its action on
$\cplane^4$ by $\mu_r(a,b,c,d)$ if the action is given by
$$
\xi\cdot(x,y,z,t)=(\xi^ax,\xi^by,\xi^cz,\xi^dt).
$$
Let $\tilde W_r\subset \cplane^4$ be the complex hypersurface
given by
$$
\tilde W_r=\{(x,y,z,t)|xy-z^{2r}+t^2=0 \}, r\geq 1.
$$
It has an isolated singularity at the origin. We call $\tilde W_r$
the {\em local $r$-conifold.}
Set
$$
\tilde W_r^\circ=\tilde W_r\setminus\{0\}.
$$

It is clear that, for any integer
 $a$ that is
prime to $r$, the action
$\mu_r(a,-a,1,0)$ preserves $\tilde W_r$. Set
$$
W_r=\tilde W_r/\mu_r,\;\; W_r^\circ=\tilde W_r^\circ/\mu_r.
$$
We call $W_r$ the {\em local $r$-orbi-conifold.}
Let $\tilde \omega^\circ_{r,w}$ be the symplectic structure
on $\tilde W_r^\circ$ induced from $\cplane^4$. It yields a
symplectic structure $\omega^\circ_{r,w}$ on $W^\circ_r$.

\subsection{The small resolutions of $W_r$ and flops}\label{sect_1.2}
By blow-ups, we have two small resolutions of $\tilde W_r$. They are
\begin{eqnarray*}
\tilde{W}^s_r &=&\{((x,y,z,t),[p,q])\in \cplane^4\times \pone^1 \\
     && |xy-z^{2r}+t^2=0,\;\;\frac{p}{q}=\frac{x}{z^r-t}=\frac{z^r+t}{y}
     \}\\
\tilde{W}^{sf}_r &=&\{((x,y,z,t),[p,q])\in \cplane^4\times \pone^1 \\
     && |xy-z^{2r}+t^2=0,\;\;\frac{p}{q}=\frac{x}{z^r+t}=\frac{z^r-t}{y}
     \}.
\end{eqnarray*}
Let
$$
\tilde \pi_r^s:\tilde W^s_r\to W^s_r, \;\;\;
\tilde \pi_r^{sf}:\tilde W^{sf}_r\to W^{sf}_r
$$
be the projections.
The extremal rays $(\tilde\pi_r^s)^{-1}(0)$
and $(\tilde\pi_r^{sf})^{-1}(0)$ are denoted by $\tilde\Gamma^s_r$ and
$\tilde\Gamma^{sf}_r$ respectively.
Both of them are isomorphic to $\pone^1$.
The action of $\mu_r$ extends naturally to
both resolutions by setting
$$
\xi\cdot[p,q]=[\xi^ap,q]
$$
for the first model and
$$
\xi\cdot[p,q]=[\xi^{-a}p,q]
$$
for the second one.

Set
$$
W^s_r=\tilde{W}^s_r/\mu_r,
\;\;\;W^{sf}_r=\tilde{W}^{sf}_r/\mu_r
\;\;\;\Gamma^s_r=\tilde{\Gamma}^s_r/\mu_r
\;\;\;\Gamma^{sf}_r=\tilde{\Gamma}^{sf}_r/\mu_r.
$$
We call $W^s$ and $W^{sf}$  {\em small resolutions} of $W_r$.
We say that $W^{sf}$ is the flop of $W^s$ and vice versa. They are
both orbifolds with singular points on $\Gamma^s$ and
$\Gamma^{sf}$.

Another important fact we use in this paper is
\begin{prop}\label{proposition_2.2.1}
For $r\geq 2$,
the normal bundle of $\tilde\Gamma^s_r$ {\em ($\tilde \Gamma^{sf}_r$)}
in $\tilde W^s_r$ {\em ($\tilde W^{sf}_r$)}
is $\mc O\oplus \mc O(-2)$.
\end{prop}
{\bf Proof. }The proof is given in \cite{Laufer}.

\subsection{Orbifold structures on $W^s$ and $W^{sf}$}\label{sect_1.25}
Let us take $W^s$. The singular points are points 0 and $\infty$
on $\Gamma^s$. In term of $[p,q]$ coordinates, they are
$$
0=[0,1];\;\; \infty=[1,0].
$$
We denote them by $\mfp$ and $\mfq$ respectively. Since $\tilde
W^s\subset \cplane^5$ near $\mfp$, the (tangent) of a uniformizing
system of $\mfp$ is given by
$$
\{ (p,x,y,z,t)|x=t=0\}.
$$
$\mu_r$ acts on this space by
$$
\xi(p,y,z)=(\xi^ap,\xi^{-a}y,\xi z).
$$
At $\mfp$, for each given $\xi^k=\exp(2\pi ik/r), 1\leq k\leq r$,
there is a corresponding twisted sector(\cite{CR1}). As a set, it
is same as $\mfp$. We denote this twisted sector by $[\mfp]_{k}$.
For each twisted sector, a degree shifting number is assigned. We
conclude that
\begin{lemma}\label{lemma_2.3.1}
For $\xi^k=\exp(2\pi ik/r), 1\leq k\leq r$, the degree shifting
$$
\iota([\mfp]_{k})= 1+ \frac{k}{r}.
$$
\end{lemma}
\n{\bf Proof. }This follows directly from the definition of degree
shifting. q.e.d.

\v Similar results hold for the singular point $\mfq$. Hence we
also have twisted sector $[\mfq]_k$ and
$$
\iota([\mfq]_{k})= 1+ \frac{k}{r}.
$$

A similar structure applies to $W^{sf}$. There are two singular points, denoted by
$\mfP,\mfQ$. The corresponding twisted sectors are $[\mfP]_k, [\mfQ]_k$. Then
$$
\iota([\mfP]_k)=\iota([\mfQ]_k)=1+\frac{k}{r}.
$$

\subsection{The deformation of $W_r$}\label{sect_1.3}
For convenience, we change coordinates:
$$
x=z_1+\sqrt{-1}z_2,\;\;
y=z_1-\sqrt{-1}z_2,\;\;z=\sqrt[2r]{-1}z_3,\;\;t=z_4.
$$
Thus in terms of the new coordinates $\tilde W_r$ is given by a new
equation
\begin{equation}\label{eqn_2.1}
z_1^2+z_2^2+z_3^{2r} +z_4^2=0.
\end{equation}
It is also convenient to use real coordinates
$$
(x_1,y_1,x_2,y_2,x_3,y_3,x_4,y_4)=(z_1,z_2,z_3,z_4).
$$
In terms of real coordinates, $\mu_r(a,-a,1,0)$ action is given by
$$
\rroot\cdot\left( \begin{array}{l} x_1\\ y_1\\x_2\\y_2\end{array}\right)
=\left ( \begin{array}{cccc}
 \cos\frac{2\pi a}{r}& 0& -\sin\frac{2\pi a}{r}& 0\\
 0& \cos\frac{2\pi a}{r}& 0& -\sin\frac{2\pi a}{r}\\
 \sin\frac{2\pi a}{r}& 0&\cos\frac{2\pi a}{r}& 0\\
0&\sin\frac{2\pi a}{r}& 0&\cos\frac{2\pi a}{r}\\
\end{array}\right)\left(\begin{array}{l} x_1\\y_1\\x_2\\y_2\\
\end{array}\right),
$$
and
$$
\rroot\cdot\left( \begin{array}{l} x_3\\ y_3\\x_4\\y_4\end{array}\right)
=\left( \begin{array}{cccc}
 \cos\frac{2\pi }{r}&  -\sin\frac{2\pi }{r}& 0&0\\
 \sin\frac{2\pi }{r}&\cos\frac{2\pi }{r}& 0&0\\
0&0& 1&0\\0&0&0&1\\
\end{array}\right)\left(\begin{array}{l} x_3\\y_3\\x_4\\y_4\\
\end{array}\right).
$$
The equation for $\tilde W_r$ is
$$
\left\{
\begin{array}{l}
x_1^2+x_2^2+f^2(x_3,y_3)+x_4^2=y_1^2+y_2^2+g^2(x_3,y_3)+y_4^2\\
x_1y_1+x_2y_2+f(x_3,y_3)g(x_3,y_3)+x_4y_4=0.
\end{array}
\right.
$$
Here $f$ and $g$ are defined by
$$
f(x,y)+\sqrt{-1}g(x,y)=(x+\sqrt{-1}y)^r.
$$
We propose
\begin{defn}\label{defn_2.4.1}
The {\em deformation} of $\tilde W_r$ is the set $\tilde Q_r$
defined by
$$\left\{ \begin{array}{l}
 x_1^2+x_2^2+f^2(x_3,y_3)+x_4^2=1,\\
x_1y_1+x_2y_2+f(x_3,y_3)g(x_3,y_3)+x_4y_4=0.
\end{array} \right. $$
The action $\mu_r(a,-a,1,0)$
preserves $\tilde Q_r$. Hence we set
$$
Q_r=\tilde Q_r/\mu_r
$$
and called it the {\em deformation of $W_r$}.
\end{defn}
\begin{lemma}\label{lemma_2.4.1}
$\tilde{Q}_r$ is a 6-dimensional symplectic
submanifold of $\rone^4\times \rone^4.$
\end{lemma}
{\bf Proof.}
Consider the map
$$\begin{array}{rl}
&F :\;\; \rone^4\times \rone^4 \rightarrow \rone^2\\
&(x,y)\to (F_1(x,y),
F_2(x,y)).\end{array}$$
given by
\begin{eqnarray*}
F_1(x,y)&=&x_1^2+x_2^2+f^2(x_3,y_3)+x_4^2-1,\\
F_2(x,y) &=&
x_1y_1+x_2y_2+f(x_3,y_3)g(x_3y_3)+x_4y_4.
\end{eqnarray*}
Then
$F^{-1}(0)=\tilde{Q}_r.$
The Jacobian  of $F$ is
$$
\left(\begin{array}{cccccccc}
2x_1&2x_2&2f\frac{\partial f}{\partial
x_3}&2x_4&0&0&2f\frac{\partial
f}{\partial y_3}&0\\
y_1&y_2&g\frac{\partial f}{\partial x_3}+f\frac{\partial
g}{\partial x_3}&y_4&x_1&x_2&g\frac{\partial f}{\partial
y_3}+f\frac{\partial g}{\partial y_3}&x_4\\\end{array}\right).$$
We claim that this is a rank 2 matrix:
if one of $x_1,x_2, x_4$, say $x_i$, is nonzero, the above matrix
has a rank 2 submatrix
$$\left (\begin{array}{ll}
2x_i&0\\y_i&x_i \end{array}\right).$$
Otherwise, say $(x_1,x_2, x_4)=(0,0,0);$
then by the definition of $\tilde{Q}_r$ we have $f(x_3,y_3)\not= 0,$
and $g(x_3,y_3)=0.$  Then since $f+\sqrt{-1}g$ is a
 holomorphic function of $x_3+\sqrt{-1}y_3$, we
have
$$\left|\begin{array}{cc}
 2f\frac{\partial f}{\partial
x_3}& 2f\frac{\partial
f}{\partial y_3} \\
 g\frac{\partial f}{\partial x_3}+f\frac{\partial
g}{\partial x_3} &g\frac{\partial f}{\partial y_3}+f\frac{\partial
g}{\partial y_3} \\\end{array}\right|=(\frac{\partial f}{\partial
x_3})^2 + (\frac{\partial f}{\partial y_3})^2 \not= 0.$$
Hence
$F$ has rank 2 everywhere on $\tilde Q_r$ and 0 is its regular value.
This     implies
that
$\tilde{Q}_r$ is a smooth 6-dimensional submanifold of $\rone^4\times \rone^4.$\\

Next we prove that
$\tilde{Q}_r$ has a canonical symplectic structure
$\omega_{\tilde{Q}_r}$ induced from
$$(\rone^4\times \rone^4,\omega_o=- \Sigma
dx_i\wedge dy_i).$$
It is sufficient to prove that
$$
\omega_o(\nabla F_1, \nabla F_2)\neq 0.
$$
By direct computations,
\begin{eqnarray*}
\nabla F_1&=&(2x_1, 2x_2,
2f\frac{\partial f}{\partial x_3}, 2x_4,0,0,2f\frac{\partial f}{\partial
y_3},y_3 ),\\
\nabla F_2&=&(y_1,y_2,f\frac{\partial g}{\partial
x_3}+g\frac{\partial f}{\partial x_3},y_4,x_1,x_2,f\frac{\partial
g}{\partial y_3}+g\frac{\partial f}{\partial y_3},x_4),
\end{eqnarray*}
Therefore
\begin{eqnarray*}
 -\omega_o(\nabla F_1, \nabla F_2)&=&\sum dx_i(\nabla F_1)dy_i(\nabla
F_2)-dx_i(\nabla F_2)dy_i(\nabla
F_1)\\
&=&2x_1^2+2x_2^2+2f((\frac{\partial f}{\partial
x_3})^2+(\frac{\partial g}{\partial x_3})^2)+2x_4^2\neq 0.
\end{eqnarray*}
  Hence $\tilde{Q}_r$ is
   a symplectic submanifold with a canonical symplectic structure
 induced from $\rone^4\times \rone^4.$
 q.e.d.

\v\n
We denote the symplectic structure by $\tilde\omega^\circ_{r,q}$.

 Put
$$
\tilde{L}_r:= \{(x,y)\in \tilde{Q}_r |y_1=y_2=g(x_3,y_3)=y_4=0 \}.
$$
and set
$$
\tilde Q_r^\circ=\tilde Q_r\setminus \tilde L_r.
$$
The $\mu_r$-action preserves $\tilde L_r$; we set
$$
L_r=\tilde L_r/\mu_r,\;\; Q_r^\circ =\tilde Q_r^\circ/\mu_r.
$$
$L_r$ is  the exceptional set in $ Q_r$ with respect to the
deformation in the following sense:
\begin{lemma}\label{lemma_2.4.2}
There is a natural diffeomorphism between $W_r^\circ$
and $Q_r^\circ$.
\end{lemma}
{\bf Proof. } We denote by
$[x,y] \in W_r^\circ$
the equivalence class of $(x,y)\in \tilde{W}_r$ with respect to the
quotient by $\mu_r$.

For any $\lambda>0$ we let $\tilde W_{r,\lambda}\subset \tilde W_r$
be the set of $(x,y)$ satisfying
$$
x_1^2+x_2^2+f^2(x_3,y_3)+x_4^2
 =y_1^2+y_2^2+g^2(x_3,y_3)+y_4^2 =\lambda
$$
and
$$
x_1y_1+x_2y_2+ f(x_3,y_3)g(x_3,x_3) +x_4y_4=0.
$$
It is not hard to see that
\begin{itemize}
\item $\tilde W_{r,\lambda}$ is preserved by the $\mu_r$ action; set
$$
W_{r,\lambda}=\tilde W_{r,\lambda}/\mu_r;
$$
\item $\tilde W_r^\circ $ is foliated by $\tilde W_{r,\lambda},\lambda
\in \rone^+$.
\end{itemize}

On the other hand, $\tilde Q_r^\circ$ has a similar foliation:
for $\lambda >0$, let $\tilde Q_{r,\lambda}\subset \tilde Q_r$ be
the set of $(x,y)$ satisfying
\begin{eqnarray*}
&& x_1^2+x_2^2+f^2(x_3,y_3)+x_4^2=1,\\
&& y_1^2+y_2^2+g^2(x_3,y_3)+y_4^2 =\lambda^2,\\
&& x_1y_1+x_2y_2+ f(x_3,y_3)g(x_3,x_3) +x_4y_4=0.
\end{eqnarray*}
Then
\begin{itemize}
\item $\tilde Q_{r,\lambda}$ is preserved by the $\mu_r$ action; set
$$
Q_{r,\lambda}=\tilde Q_{r,\lambda}/\mu_r;
$$
\item $\tilde Q_r^\circ $ is foliated by $\tilde Q_{r,\lambda},\lambda
\in \rone^+$.
\end{itemize}
We next introduce the identification between $W_{r,\lambda}$
and $Q_{r,\lambda}$. Let
$u_\lambda(x_3,y_3)$ and $ v_\lambda(x_3,y_3)$ be functions that solve
$$
(u+iv)^r=\lambda^{-1}f(x_3,y_3)+\sqrt{-1} \lambda g(x_3,y_3).
$$
Such a pair $u+iv$ exists up to a factor $\xi^k$. Then
\begin{eqnarray*}
&[x_1,x_2,x_3,x_4,  y_1,y_2,y_3,y_4]
\longleftrightarrow \\
&[\lambda^{-1}x_1,\lambda^{-1}x_2, u(x_3,y_3), \lambda^{-1}x_4,
\lambda y_1,\lambda y_2,v(x_3,y_3),\lambda y_4]
\end{eqnarray*}
induces an identification between $W_{r,\lambda}$ and $Q_{r,\lambda}$,
and therefore between $W_r^\circ$ and $Q_r^\circ$. q.e.d.

\v\n We denote the identification map constructed in the proof by
$$
\Phi_r: W_r^\circ \to Q_r^\circ.
$$
In particular, we note that
the restriction of $
\Phi_r$ to ${W_{r,1}}$ is the identity.

\subsection{The comparison between local $r$-orbi-conifolds and local conifolds}\label{sect_1.4}
When $r=1$, the local model  is the well-known conifold.
Since $\mu_r=\mu_1=\{1\}$ is trivial, there is no orbifold structure.
It is well known that
\begin{itemize}
\item
$W_1^s$ and $W_1^{sf}$ are
$$
\mc O(-1)\oplus \mc O(-1)\to \pone^1,
$$
where $\Gamma^s$ and $\Gamma^{sf}$ are the zero section $\pone^1$;
They are flops of each other;
\item
$Q_1$ is diffeomorphic to the cotangent bundle of $S^3$. The induced
symplectic structure from $\rone^4\times \rone^4$ coincides with
the canonical symplectic structure on $T^\ast S^3$.
\item the map
$$
\Phi_1:(W_1,\omega^\circ_{1,w})\to (Q_1,\omega^\circ_{1,q})
$$ is a symplectomorphism.
\end{itemize}

There are  natural (projection) maps
$$
\pi_{r,w}: \tilde W_r\to  W_1, \;\;
\pi_{r,q}: \tilde Q_r\to Q_1
$$
given by
$$
x_i\to x_i, \;\; y_i\to y_i, i\not=3,
$$
and
$$
(x_3,y_3)\to (f(x_3,y_3),g(x_3,y_3)).
$$
Similarly, there are maps
$$
\pi_{r,w}^s:\tilde W^s_r\to W^s_1,\;\;
\pi_{r,w}^{sf}:\tilde W^{sf}_r\to W^{sf}_1.
$$
We note that all these projection maps are almost
$r$-coverings. They are coverings except on $x_3=y_3=0$, where the maps are
only $r$-branched coverings. Note  that
$$
\tilde L_r=\pi_{r,q}^{-1} L_1.
$$
It is the union of $r$ copies of $S^3$ intersecting
at
$$
\left\{
\begin{array}{l}
x_1^2+x_2^2+x_4^2=1\\
x_1y_1+x_2y_2+x_4y_4=0
\end{array}\right\}
\cap \{x_3=y_3=0\}.
$$

\section{Cohomologies  }\label{sect_2}

\subsection{Definitions}\label{sect_2.1}

 Let $(\Omega^\ast(\tilde W^\circ_r),d)$ be the
de Rham complex of $\tilde W^\circ_r$. $\mu_r$ has a natural
representation on this complex.
let
$$
\Omega^\ast_{\mu_r}(\tilde W^\circ_r)
\subset \Omega^\ast(\tilde W^\circ_r)
$$
be the subcomplex of $\mu_r$-invariant forms. We have
$$
H^\ast(W^\circ_r)=H^\ast(\Omega^\ast_{\mu_r} (\tilde
W^\circ_r),d ).
$$
Similar definitions  apply to
$W^s_r,W^{sf}_r, Q_r^\circ, Q_r, W_{r,1}=Q_{r,1}$  etc.

Then
\begin{lemma}\label{lemma_3.1.1}
$H^\ast(W^\circ_r)=H^\ast(W_{r,1})$.
\end{lemma}
{\bf Proof. }
We note that there is a $\mu_r$-isomorphism
$$
\tilde W_{r}^\circ\cong \tilde W_{r,1}\times \rone^+.
$$
In fact,  it is induced by a natural identification
\begin{eqnarray*}
\tilde W_{r,\lambda} &\leftrightarrow&\tilde W_{r,1}\times
\{\lambda\};\\
x_i&\leftrightarrow& \lambda^{-\frac{1}{2}}x_i, i\not=3;\;
x_3\leftrightarrow \lambda^{-\frac{1}{2r}}x_3,\\
y_i&\leftrightarrow& \lambda^{-\frac{1}{2}}y_i, i\not= 3;\;
y_3\leftrightarrow \lambda^{-\frac{1}{2r}}y_3.
\end{eqnarray*}
Hence $\tilde W^{\circ}_r$ is $\mu_r$-homotopy equivalent to
$\tilde W_{r,1}$. Hence the claim follows.  q.e.d.

\v\n
The result also follows from
$$
W_r^\circ\cong W_{r,1}\times \rone^+
$$
directly.
Similarly, we have
$$
Q_r^\circ \cong Q_{r,1}\times \rone^+.
$$
Hence
$$
H^\ast(Q^\circ_r)=H^\ast(Q_{r,1}).
$$
Note that $Q_{r,1}=W_{r,1}$. We have
$$
H^\ast(W^\circ_r)=H^\ast(W_{r,1})=H^\ast(Q_{r,1})=H^\ast(Q^\circ_r).
$$

\subsection{Computation of cohomologies}\label{sect_2.2}
We first study $H^\ast(W_{r,1})$.

Recall that we have a map
$$
\pi_{r,w}:\tilde W_{r,1}\to W_{1,1}
$$
given by
$$
\pi_{r,w}(x,y)=
(x_1,x_2,f(x_3,y_3),x_4,y_1,y_2,g(x_3,y_3),y_4).
$$
We now introduce a $\mu_r$ action on $W_{1,1}$. For convenience,
we use coordinates $(u,v)$ for the $\rone^4\times \rone^4$ in which
$W_{1,1}$ is embedded. Then
$$
\rroot\cdot\left( \begin{array}{l} u_1\\ v_1\\u_2\\v_2\end{array}\right)
=\left ( \begin{array}{cccc}
 \cos\frac{2\pi a}{r}& 0& -\sin\frac{2\pi a}{r}& 0\\
 0& \cos\frac{2\pi a}{r}& 0& -\sin\frac{2\pi a}{r}\\
 \sin\frac{2\pi a}{r}& 0&\cos\frac{2\pi a}{r}& 0\\
0&\sin\frac{2\pi a}{r}& 0&\cos\frac{2\pi a}{r}\\
\end{array}\right)\left(\begin{array}{l} u_1\\v_1\\u_2\\v_2\\
\end{array}\right),
$$
and acts trivially on $u_3,v_3, u_4$ and $v_4$.
Then it is clear that $\pi_{r,w}$ is $\mu_r$-equivariant.
It induces a morphism between complexes
\begin{equation}\label{eqn_3.1}
\pi_{r,w}^\ast: (\Omega^\ast_{\mu_r}(W_{1,1}),d)
\to (\Omega^\ast_{\mu_r}(\tilde W_{r,1}),d).
\end{equation}
Here $\Omega_G$ always represents the subspace that is $G$-invariant
if $\Omega$ is a $G$-representation.

\begin{prop}\label{prop_3.2.1}
$\pi_{r,w}^\ast$ in \eqref{eqn_3.1} is an isomorphism between
the {\em cohomologies} of the two complexes.
\end{prop}
{\bf Proof. }The idea of the proof is to consider a larger
connected Lie group action on spaces: Let $S^1=\{e^{2\pi i\theta}\}$.
Suppose its action on $(x,y)$ is given by
$$
e^{2\pi i\theta}\cdot\left( \begin{array}{l} x_1\\ y_1\\x_2\\y_2\end{array}\right)
=\left ( \begin{array}{cccc}
 \cos \theta& 0& -\sin\theta& 0\\
 0& \cos\theta& 0& -\sin\theta\\
 \sin\theta& 0&\cos\theta& 0\\
0&\sin\theta& 0&\cos\theta\\
\end{array}\right)\left(\begin{array}{l} x_1\\y_1\\x_2\\y_2\\
\end{array}\right),
$$
and the trivial action on $x_3,y_3,x_4$ and $y_4$.
The same action is defined on $(u,v)$. Again, $\pi_{r,w}$ is
$S^1$-equivariant.

Since $S^1$ is a connected Lie group and its actions commutes with
$\mu_r$-actions on both spaces, the subcomplex
$$
((\Omega^\ast_{\mu_r}(\tilde W_{r,1}))_{S^1},d)\subset
(\Omega^\ast_{\mu_r}(\tilde W_{r,1}),d)
$$
of $S^1$-invariant forms yields same cohomology as the original one, i.e,
$$
H^\ast((\Omega^\ast_{\mu_r}(\tilde W_{r,1}))_{S^1},d)
=H^\ast (\Omega^\ast_{\mu_r}(\tilde W_{r,1}),d)
$$
Similarly,
$$
H^\ast((\Omega^\ast_{\mu_r}(W_{1,1}))_{S^1},d)
=H^\ast (\Omega^\ast_{\mu_r}(W_{1,1}),d)
$$
It is then sufficient to show that
\begin{equation}\label{eqn_3.2}
\pi_{r,w}^\ast:H^\ast((\Omega^\ast_{\mu_r}(W_{1,1}))_{S^1},d)
 \to H^\ast((\Omega^\ast_{\mu_r}(\tilde W_{r,1}))_{S^1},d)
\end{equation}
is an isomorphism. By the definition of the actions, we note that
\begin{equation}\label{eqn_3.3}
(\Omega^\ast_{\mu_r}( W_{1,1}))_{S^1}
=\Omega^\ast_{S^1}(W_{1,1}).
\end{equation}

We now show \eqref{eqn_3.2}. Recall that $\pi_{r,w}$ is
an $r$-branched covering ramified over
$$
R_1=\left\{\begin{array}{l}
u_1^2+u_2^2+u_4^2=v_1^2+v_2^2+v_4^2=1\\
u_1v_1+u_2v_2+u_4v_4=0
\end{array}
\right\}\cap \{u_3=v_3=0\}
$$
Set $\tilde R_r=\pi^{-1}_{r,w}(R_1)$ and
$$
\tilde U_r= \tilde W_{r,1}\setminus\tilde R_r,\;\;\;
U_1=W_{1,1}\setminus R_1.
$$
  Then
$\pi_{r,w}:\tilde  R_r\to R_1$
is 1-1 and $\pi_{r,w}: \tilde U_r\to U_1$
is  an$r$-covering.

Let $V_1$ be an $S^1$-invariant tubular neighborhood of $R_1$ in $W_{1,1}$.
By the implicit function theorem, we know that
$$
V_1\cong R_1\times D_1,
$$
where $D_1$ is the unit disk in the complex plane $\cplane=\{u_3+\sqrt{-1}v_3\}$.
Let $\tilde V_r=\pi_{r,w}^{-1}(V_1)$. Then
$$
\tilde V_r\cong \tilde R_r\times D_1,
$$
where $D_1$ is the unit disk in the complex plane $\cplane=\{x_3+iy_3\}$.
In terms of these identifications, $\pi_{r,w}$ can be rewritten as
\begin{eqnarray*}
&&\pi_{r,w}: \tilde R_r\times D_1\to R_1\times D_1\\
&&\pi_{r,w}(\gamma,z_3)=(\gamma,z_3^r),
\end{eqnarray*}
where $\gamma\in \tilde R_r
=R_1, z_3=x_3+iy_3 .$

Consider the short exact sequences
$$
0 \to (\Omega^\ast_{\mu_r}( W_{1,1}))_{S^1}
\to (\Omega^\ast_{\mu_r}(U_1))_{S^1}\oplus
 (\Omega^\ast_{\mu_r}(V_1))_{S^1}
\to (\Omega^\ast_{\mu_r}(U_1\cap V_1))_{S^1}
\to 0
$$
and
$$
0 \to (\Omega^\ast_{\mu_r}(\tilde W_{r,1}))_{S^1}
\to (\Omega^\ast_{\mu_r}(\tilde U_r))_{S^1}\oplus
 (\Omega^\ast_{\mu_r}(\tilde V_r))_{S^1}
\to (\Omega^\ast_{\mu_r}(\tilde U_r\cap \tilde V_r))_{S^1}
\to 0.
$$
$\pi_{r,w}^\ast$ is a morphism between two complexes. We assert that
\begin{eqnarray}
&&\pi_{r,w}^\ast: H^\ast((\Omega^\ast_{\mu_r}(U_1))_{S^1},d)
\xrightarrow{\cong} H^\ast((\Omega^\ast_{\mu_r}(\tilde U_r))_{S^1},d),\label{eqn_3.4}\\
&&\pi_{r,w}^\ast: H^\ast((\Omega^\ast_{\mu_r}(V_1))_{S^1},d)
\xrightarrow{\cong} H^\ast((\Omega^\ast_{\mu_r}(\tilde V_r))_{S^1},d), \label{eqn_3.5}\\
&&\pi_{r,w}^\ast: H^\ast((\Omega^\ast_{\mu_r}(U_1\cap V_1))_{S^1},d)
\xrightarrow{\cong} H^\ast((\Omega^\ast_{\mu_r}(\tilde U_r\cap \tilde V_1))_{S^1},d).
\label{eqn_3.6}
\end{eqnarray}
Once these are proved, by the five-lemma, we know that
$$
\pi_{r,w}^\ast: H^\ast((\Omega^\ast_{\mu_r}(W_{1,1}))_{S^1},d)
\xrightarrow{\cong} H^\ast((\Omega^\ast_{\mu_r}(\tilde W_{r,1}))_{S^1},d)
$$
which is \eqref{eqn_3.2}.

We now explain \eqref{eqn_3.4}, \eqref{eqn_3.5} and \eqref{eqn_3.6}.
\v\n
{\em The proof of \eqref{eqn_3.4}. } We observe that
$$
\pi_{r,w}^\ast:(\Omega^\ast_{\mu_r}(U_1))_{S^1}
\xrightarrow{\cong} (\Omega^\ast_{\mu_r}(\tilde U_r))_{S^1}.
$$
Hence it induces an isomorphism on cohomology level.

\v\n
{\em The proof of \eqref{eqn_3.5}. }Since $\tilde V_r$ is
$\mu_r\times S^1$-homotopy equivalent to $\tilde R_r$,
we have
$$
H^\ast((\Omega^\ast_{\mu_r}(\tilde V_r))_{S^1},d)\cong
H^\ast((\Omega^\ast_{\mu_r}(\tilde R_r))_{S^1},d).
$$
Similarly,
$$
H^\ast((\Omega^\ast_{\mu_r}(V_1))_{S^1},d)\cong
H^\ast((\Omega^\ast_{\mu_r}(R_1))_{S^1},d).
$$
Because
$$
H^\ast((\Omega^\ast_{\mu_r}(\tilde R_r))_{S^1},d)
=H^\ast((\Omega^\ast_{\mu_r}(R_1))_{S^1},d),
$$
we have \eqref{eqn_3.5}.
\v\n
{\em The proof of \eqref{eqn_3.6}. }The proof is the same as that of
\eqref{eqn_3.4}.
\v \n This completes the proof of the theorem. q.e.d.

\v
So far, we have shown that
$$
H^\ast(W_{r,1})=H^\ast (\Omega^\ast_{\mu_r}(\tilde
W_{r,1}),d)\cong H^\ast (\Omega^\ast_{\mu_r}(W_{1,1}),d)
=H^\ast((\Omega^\ast_{\mu_r}(W_{1,1}))_{S^1},d).
$$
Furthermore, by \eqref{eqn_3.3} we have
$$
H^\ast((\Omega^\ast_{\mu_r}(W_{1,1}))_{S^1},d)
=H^\ast(\Omega^\ast_{S^1}(W_{1,1}),d)= H^\ast(W_{1,1}).
$$
Since $W_{1,1}\cong S^3\times S^2$ we have
\begin{corollary}\label{cor_3.2.1}
$H^\ast(W_{r,1})\cong H^\ast(S^3\times S^2)$.
\end{corollary}

Let $H_1$ be a generator of $H^2(S^3\times S^2)$
such that
$$
\int_{S^2} H_1=1.
$$
Here $S^2$ is any fiber  $\{x\}\times S^2$ in
$S^3\times S^2$. Set
$$
\tilde H_r=\pi_{r,w}^\ast H_1
$$
and let $H_r$ be its induced form on $W_{r,1}$. This is a generator of
$H^2(W_{r,1})$. Without loss of generality, we also assume that
it is a generator of $H^2(W_{r}^\circ)$.

Let $\omega_{r,w}$ and $\omega_{r,q}$ be symplectic forms
on $W_r^\circ$ and $Q_r^\circ$ respectively. Suppose that
$$
[\omega_{r,w}|_{W_{r,1}}]=[\omega_{r,q}|_{Q_{r,1}}].
$$
Here $[\omega]$ denotes the cohomology class of $\omega$. Then
there exists a symplectomorphism
$$
\Phi'_r:( W_r^\circ, \omega_{r,w})\to (Q_r^\circ,\omega_{r,q}).
$$
In fact, by the assumption, we have
$$
[\omega_{r,w}]=[\Phi_r^\ast\omega_{r,q}].
$$
Then, by the standard Moser argument, there exists a diffeomorphism
$$
f: W_r^\circ\to W_r^\circ
$$
such that $f^\ast \omega_{r,w}= \Phi_r^\ast\omega_{r,q}$. Now we can
set $\Phi'_r=\Phi_r\circ f^{-1}$.  In particular, by applying it
to $\omega_{r,w}^\circ$ and $\omega_{r,q}^\circ$ we have
\begin{corollary}\label{cor_3.2.2}
There exists a symplectomorphism
$$
\Phi'_r:( W_r^\circ, \omega^\circ_{r,w})\to (Q_r^\circ,\omega^\circ
_{r,q}).
$$
\end{corollary}
{\bf Proof. }We observe that both  symplectic forms are exact.
Hence they represent the same cohomology class, namely 0. q.e.d.

\v\n Next we consider $H^\ast(W^s_r)$. The argument is same
as above: we also have a map
$$
\pi_{r,w}: \tilde W^s_r\to W^s_1.
$$
This map will induce an isomorphism
\begin{prop}\label{prop_3.2.2}
$H^\ast(W^s_r)=H^\ast(W^s_1)$.
\end{prop}
{\bf Proof. }Since the proof is parallel to that of
proposition \ref{prop_3.2.1}, we only sketch the proof.

We use complex coordinates $(x,y,z,t,[p,q])$ for
$\tilde W^s_r$ and $(u,v,w,s,[m,n])$ for $W^s_1$.
Then $\pi_{r,w}$ is induced by the map
$$
u=x,\; v=y,\; w=z^r,\; s=t,\; \frac{m}{n}=\frac{p}{q}.
$$
We can introduce a $\mu_r$-action on $W^s_1$ by
$$
\xi(u,v,w,s,[m,n])=
(\xi^a u,\xi^{-a}v,w,s,[\xi^am,n]), \xi=\rroot.
$$
Then $\pi_{r,w}$ is $\mu_r$-equivariant.

Moreover, both spaces admit an $S^1$-action such that
$\pi_{r,w}$is $S^1$-equivariant: for $\xi\in S^1$:
\begin{eqnarray*}
\xi(x,y,z,t,[p,q])&=& (\xi^ax,\xi^{-a}y, z,t,[\xi^ap,q])\\
\xi(u,v,w,s,[m,n])&=&
(\xi^a u,\xi^{-a}v,w,s,[\xi^am,n]).
\end{eqnarray*}
$\pi_{r,w}$ is an $r$-branched covering ramified over
$$
W^s_1\cap \{w=0\}.
$$
Then the rest of the proof is simply a copy of   the argument in
Proposition \ref{prop_3.2.1}.
 q.e.d.

\v
Since
$$
W^s_1\cong \mc O(-1)\oplus \mc O(-1),
$$
$H^2(W^s_1)=H^2(\pone^1)$ is 1-dimensional. So is
$H^2(W^s_r)$. Let $H^s_{r}$ be the generator of $H^2(W^s_r)$
such that
$$
\int_{\Gamma^s_r}H^s_r=1.
$$
Since the normal bundle of $\tilde \Gamma_{r}^s$ is $\mc O\oplus
\mc O(-2)$, it admits a symplectic form $\omega'$. We normalize
it by
$$
\int_{\Gamma^s_r}\omega'=1.
$$
It induces a symplectic structure, denoted by
$\omega^s_{r}$ on the neighborhood $U$ of $\Gamma^s_r$.
It is easy to see that
this symplectic structure is tamed by the complex structure on
$U$. Hence we conclude that
\begin{corollary}\label{cor_3.2.3}
There is a symplectic form on $W^s_r$
that represents the class $H^s_r$
and is tamed by its complex structure. This form is denoted by
$\omega^s_r$.
\end{corollary}

\section{Orbifold symplectic flops}\label{sect_3}

\subsection{The global orbi-conifolds}\label{sect_3.1}

Following  \cite{STY} we give the definition of orbi-conifolds.
 \begin{defn}\label{defn_4.1.1}
A {\em real 6-dimensional orbi-conifold }
is a topological space $Z$
covered by an atlas of charts $\{(U_i,\phi_i)\}$ of the following
two types: either $(U_i,\phi_i) $ is an  orbifold chart or
$$\phi_j: U_j \rightarrow W_{r_j}$$
is a
 homeomorphism onto   $W_{r_j}$ defined in  \S\ref{sect_1.1}.
 In the latter case, we call the point $\phi_j^{-1}(0)$ a
 singularity of $Z$.

 Moreover, the transition maps
$\phi_{ij}=\phi_i\circ \phi_j^{-1}$ must be
 smooth in the orbifold sense  away from singularities
 and if $p\in U_i \cap U_j$ is a singularity then we have $r_i=r_j$ (denote
it by $r$), and
 there must be
 an open subset $N\subset \cplane^4$ containing 0 such that the lifting
of $\phi_{ij}, $
 $$
\tilde{\phi}_{ij}:
\tilde{W}_{r}\cap N\longrightarrow \tilde{W}_{r}\cap N
$$
 in the uniformizing system is
 the restriction of an analytic isomorphism
$\tilde{\phi}:\;\cplane^4 \rightarrow \cplane^4$
which is smooth   away from the origin, $C^1$ at the origin with
 $d\tilde{\phi}_0 \in Sp(8,\rone)$, and set-wise fixes $\tilde{W}_r$.
\end{defn}
 We call such charts {\it smooth admissible coordinates.} Note that
 in the case $r=1$ the singularity is the ordinary double point discussed in
 [STY].

>From now on, we label the set of singularities
$$
P=\{p_1,p_2,\ldots\},
$$
and for each point $p_i$ its local model is given by a
standard
model $W_{r_i}$.

\begin{defn}\label{defn_4.1.2}
A {\em symplectic
structure} on an orbi-conifold $Z$ is a smooth orbifold symplectic
form $\omega_Z$ on the orbifold $Z\setminus P$
which,   around each
singularity $p_i$, coincides with  $\omega^\circ_{w,r_i}$.
$(Z,\omega_Z)$ is called a {\em symplectic orbi-conifold}.
\end{defn}

>From now on, we assume that $Z$ is compact and
$|P|=\kappa$. One can perform a smoothing   for each
 singularity  of $Z$ as in \S\ref{sect_1.3}
- replace a neighborhood  of each singularity
 $p_i$ by  a neighborhood of $L_{r_i}$ in $Q_{r_i}$ - to get an orbifold. We denote
this orbifold  by $X$.


For each singularity $p_i$ of $Z$ we can perform two
  small resolutions, i.e., we replace the neighborhood of the
singularity by $W_{r_i}^s$
  or $W_{r_i}^{sf}$ as in \S\ref{sect_1.2}.
  There are
  $2^\kappa$ choices of small resolutions, and so we get
  $2^\kappa$ orbifolds $ Y_1, \cdots, Y_{2^\kappa}.$
\begin{defn}\label{defn_4.1.3}
Two small resolutions $Y$ and $Y'$ are said to be {\em flops of
each other} if at each $p_i$, one is obtained by replacing
$W_{r_i}^s$ and the other by $W_{r_i}^{sf}$. We write
$Y'=Y^f$ and vice versa.
\end{defn}

\subsection{Symplectic structures on $Y_i$'s and flops}\label{sect_3.2}

Not every small resolution $Y_\alpha, 1\leq \alpha\leq 2^\kappa$
admits a symplectic structure. Our first main theorem of the paper   gives a necessary and sufficient
condition for $Y$ to have a symplectic structure in terms
of the topology of $X$.

Let $L_{r_i}\subset X$. For simplicity, we assume its neighborhood to be
$Q_{r_i}$. Recall that there is a projection map
$$
\pi_{r_i,q}: \tilde Q_{r_i}\to Q_1.
$$
Let $\Theta_1$ be the Thom form of the normal bundle of
$L_1=S^3$ in $Q_1$.
We assume it is supported in a small neighborhood
of $L_1$. Set
$$
\tilde\Theta_{r_i}=\pi_{r_i,q}^\ast \Theta_1.
$$
We can choose $\Theta_1$ properly such that $\tilde \Theta_{r_i}$
is $\mu_{r_i}$ invariant. Hence it induces a local form $\Theta_{r_i}$
on $Q_{r_i}$ and hence on $X$.

Then we restate Theorem \ref{theorem_1.1}:
{\em
One of the $2^\kappa$
small resolutions
admits a symplectic stucture if and only if  on
 $X$ we have  the following  cohomology relation
\begin{equation}\label{eqn_4.1}
[\sum_{i=1}^\kappa
\lambda_i\Theta_{r_i}]=0\in  H^3(X,\rone)
\;\;with \;\lambda_i\not= 0\; for \;all\;i.
\end{equation}
}

As a corollary,
\begin{corollary}\label{cor_4.2.1}
Suppose we have a pair of resolution $Y$ and $Y^f$ that are flops
of each other. Then $Y$ admits a symplectic structure if and only
if $Y^f$ does.
\end{corollary}
$Y^f$ is then called the {\em symplectic flop } of $Y$.

\section{Proof of theorem \ref{theorem_1.1}}\label{sect_4}

\subsection{Necessity}\label{sect_4.1}
We first prove that \eqref{eqn_4.1} is necessary.

Suppose that we have a $Y$ that admits a symplectic structure
$\omega$.
For simplicity, we assume that at each singular point $p_i\in Z$,
it is replaced by $W_{r_i}^s$ to get $Y$. The extremal ray is
$\Gamma^s_i$. Set
$$
\lambda_i= \int_{\Gamma^s_{r_i}}\omega=\frac{1}{r_i}\int_{\tilde
\Gamma^s_{r_i}}\tilde \omega.
$$

Now we consider the pair of spaces $(X,X\setminus \cup L_{r_i})$. The exact
sequence of the (orbifold) de Rham complex of the pair is
$$
0 \to \Omega^{\ast-1}(X\setminus \cup L_{r_i})
\xrightarrow{\gamma} \Omega^\ast (X,X\setminus \cup L_{r_i})
\xrightarrow{\delta}
\Omega^\ast(X)\to 0.
$$
given by
$$
\gamma(f)=(0,f), \;\; \delta(\alpha,f)= \alpha.
$$
It induces a long exact sequence on (orbifold) cohomology
$$
\cdots\to
H^2(X\setminus\cup L_{r_i})\rightarrow
H^3(X,X\setminus\cup
    L_{r_i})
    \rightarrow H^3(X)\to\cdots
    $$
And applying this to $\omega$ on
$
Z\setminus P\cong X\setminus \cup_i L_{r_i},
$
we have
$$
\omega\mapsto (0,\omega) \mapsto 0.
$$
This says that
$$
[\delta\circ\gamma(\omega)]=0.
$$
We compute the left hand side of the equation.
First, by applying the excision principle
 we get
$$
H^3(X,X\setminus\cup_i L_{r_i})
\cong\bigoplus_i H^3(Q_{r_i}, Q_{r_i}^\circ) .
$$
This reduces the computation to the local case.

Let $\omega_{r_i,w}$ be the restriction of
$\omega$ in the neighborhood, simply denoted
by $W_{r_i}^s$, of $\Gamma_{r_i}^s$.
It induces a form $\omega_{r_i,q}$ on
$Q_{r_i}^\circ$. Suppose that
$$
\omega_{r_i,q}=c_i H_{r_i},
$$
where $H_{r_i}$ is the generator on $Q_{r_i,1}$, hence on
$Q_{r_i}^\circ$.
Let $\beta$ be a cut-off function such that
$$
\beta(t)=\left\{
\begin{array}{lll}
1, &  \;if \; & t>0.5;\\
0, & \;if \; & t<0.25.
\end{array}
\right.
$$
By direct computation, we have
$$
\delta\circ\gamma(H_{r_i})= d(\beta(\lambda)H_{r_i})
=\Theta_{r_i}.
$$
Therefore, we conclude that
$$
\sum_{i=1}^\kappa c_i\Theta_{r_i}=0
$$
In order to show \eqref{eqn_4.1}, it  remains to
prove that
\begin{prop}\label{prop_5.1.1}
$c_i=-\lambda_i$.
\end{prop}
{\bf Proof. }The computation is done on
$\tilde W^{s}_{r_i}$.

Take an $S^2$ in $Q_{1,1}$ as
$$
B_1=
\{(1,0,0,0,0,v_2,v_3,v_4)\in \tilde{Q}_{r_i}| v_2^2+v_3^2+v_4^2=1\}
$$
Let $\tilde B_r=\pi_{r,q}^{-1}(B_1)$.
It is
$$
\tilde B_{r_i}=
\{(1,0,x_3,0,0,y_2,y_3,y_4)\in \tilde{Q}_{r_i}| y_2^2+g^2(x_3,y_3)+v_4^2=1,
f(x_3,y_3)=0\}
$$
Then
$$
\int_{\tilde B_{r_i}}\tilde H_{r_i}=r_i\int_{B_1}H_1=r_i.
$$
Hence
$$
\int_{\tilde B_{r_i}}\omega_{r_i,q}=c_ir_i.
$$
Next we explain that
\begin{equation}\label{eqn_5.1}
\int_{\tilde B_{r_i}}\omega_{r_i,q}=-\lambda_ir_i.
\end{equation}
Then the claim follows from these two identities.

\v\n
{\em Proof of \eqref{eqn_5.1}: } We treat $B_1$ and
$\tilde B_{r_i}$ as subsets of $W^s_1$ and $\tilde W_{r_i}^s$.
By Proposition \ref{prop_3.2.1}, we assume
$\omega_{r_i,w}$ is homologous to $\pi_{r_i,w}^\ast \omega$
for some $\omega\in H^2(W^s_1)$. Then
$$
\int_{\tilde B_{r_i}}\omega_{r_i,q}=r_i
\int_{ B_{1}}\omega.
$$
On the other hand, $B_1$ is homotopic to $-\Gamma^s_1$:
via complex coordinates $W_1$ is given by
$$
uv-(w-s)(w+s)=0.
$$
The equation of the small resolution  ${W}^s_{1}$ in the
chart $\{q\not=0\}$ is
$$ \zeta v-(w-s)=0,,$$
where $  \zeta=\frac{m}{n}=\frac{u}{w+s}$ is the coordinate
of the exceptional curve $\Gamma_1^s.$
Recall that
on $B_1$ the complex coordinates are
$$
x= 1+ y_2,\;\; y= 1- y_2,\;\;  z= \sqrt{-1}y_3,\;\;
t=y_4.
$$
We have a projection map
$$B_1 \longrightarrow
\Gamma_1^s$$
given by
$$\eta =\frac{x}{z+t}=\frac{1+\sqrt{1-y_3^2-y_4^2}}
{\sqrt{-1} y_3+y_4}.$$
Here we take $y_3,\;y_4$ as coordinates on $B_1.$ It
is easy to see that this is a one to one map and the point with
$\sqrt{-1} y_3+y_4=0$ corresponds to the point "$\infty$" of
$-\Gamma_1^s.$ The sign is due to the orientation.

Let
$$
(\zeta,y,z,t)=(\frac{1+\sqrt{1-y_3^2-y_4^2}}{\sqrt{-1}y_3+y_4},
1-y_2,
   iy_3,y_4)
$$
be   any point
  in   $B_1$;  then
$$
(\zeta_0,0,0,0) =
   (\frac{1+\sqrt{1-y_3^2-y_4^2}}{\sqrt{-1}y_3+y_4},0,0,0)
$$ is in
$\Gamma_1^s.$  We construct a subset $\Lambda_1$ of $W^s_1$
$$
\rho(y_3, y_4,s)=\{(\frac{1+\sqrt{1-y_3^2-y_4^2}}{\sqrt{-1}y_3+y_4},
s(1-y_2),
  s\sqrt{-1}y_3,sy_4)\}$$
   where $0\leq s \leq1$
and $y_3,\;y_4$ are the coordinates of $N_1.$
This is a smooth 3-dimensional submanifold with
   boundary
 $$
\{\rho(y_3,y_4,0)=-\Gamma_1^s\} \cup
\{\rho(y_3,y_4, 1)=B_1\}.
$$

 It gives us a homotopy between $-\Gamma_1^s$
and $B_1$. Then
$$
\int_{\tilde B_{r_i}}{\omega_{r_i,w}}=
r_i\int_{B_1}\omega=-r_i\int_{\Gamma^s_1}\omega
=-\int_{\tilde \Gamma_{r_i}^s} \omega_{r_i,w}=-r_i\lambda_i.
$$
This shows \eqref{eqn_5.1}.
\v\n
We have completed the proof of  the proposition. q.e.d.

\v\n
This completes the proof of necessity.
\begin{remark}\label{rmk_5.1.1}
If the local resolution is $W_{r_i}^{sf}$,
$$
[\delta\circ\gamma (\omega_{r_i,w})] =\lambda_i \Theta_{r_i}.
$$
\end{remark}

\subsection{Sufficiency}\label{sect_4.2}

Suppose that \eqref{eqn_4.1} holds for  $X$: i.e, there exists
$\lambda_i$ such that
$$
\sum_i\lambda_i\Theta_{r_i}=0.
$$
For the moment we assume that all $\lambda_i<0$.
Let
$Y$ be a small resolution of $Z$ obtained by replacing
the neighborhood of $p_i$ by $W^s_{r_i}$. We assert that
$Y$ admits a symplectic structure.

>From the exact
sequence of the  pair of spaces $(X,X\setminus\cup_i L_{r_i})$
$$
H^2(X\setminus\cup_iL_{r_i})\xrightarrow{\gamma}  H^3(X,X\setminus  L_{r_i})
    \rightarrow H^3(X)
$$
we conclude that there exists a 2-form $\sigma^*\in
    H^2(X\setminus\cup_iL_{r_i})$ such that
$$
\gamma(\sigma^\ast)=\sum \lambda_i\Theta_{r_i}.
$$
 since
$$
X\setminus\cup_iL_{r_i}\cong
    Y\setminus\cup_i\Gamma_{r_i}^s,
$$
$\sigma^*\in
    H^2(Y\setminus\cup_i\Gamma_{r_i}^s).$
On the other hand, we consider the exact
sequence of the  pair of spaces $(Y,Y\setminus\cup_i \Gamma_{r_i}^s )$
$$
H^2(Y)\rightarrow H^2(Y\setminus\cup_i\Gamma_{r_i}^s)
\rightarrow  H^3(Y,Y\setminus \cup \Gamma_{r_i}^s )
\cong \bigoplus_iH^3(W_{r_i}^s, W_{r_i}^\circ ).
$$
It is known that locally $\tilde W_{r_i}^s$
is diffeomorphic to its normal bundle $\mc O\bigoplus \mc O(-2)$
     of $\tilde\Gamma_{r_i},$ thus
$$
H^3(Y,Y\setminus \cup \Gamma_{r_i}^s )=0.
$$
It follows that there exist a 2-form $\sigma\in
    H^2(Y )$ which extends  $\sigma^* $.

Let $U_i$ be a small neighborhood of $\Gamma_{r_i}^s$ in
$Y$ and $\tilde U_i\subset
\tilde W^s_{r_i}$ be its pre-image in the uniformizing
system. Set
$$
\sigma_i=\sigma|_{U_i}.
$$
By the proof of necessity, we know that
$$
[\sigma_i]= [-\lambda_i\omega^s_{r_i}].
$$
Then we can deform
$\sigma_i$ in its cohomology class near $\tilde\Gamma^s_{r_i}$
  such that
$$
\sigma_i=-\lambda_i\omega^s_{r_i}.
$$
Hence we get a new form $\sigma$ on $Y$ that gives symplectic
forms near $\Gamma_i^s$.
On the other hand, we have a  form $\omega_Z$ on $Z$ that is symplectic
away from $P$. This form extends to $Y$, still denoted
by $\omega_Z$, but is degenerate at the $\Gamma_{r_i}^s$.
For sufficiently large $N$ we have
$$
\Omega= \sigma + N\omega_Z.
$$
This is a symplectic structure on $Y$: $\Omega$ is non-degenerate
away from a small neighborhood of the $\Gamma_{r_i}^s$ for large $N$;
both $\sigma$ and $\omega_Z$ are tamed by the complex structure
in the $U_i$, i.e,
$$
\sigma(\cdot, J\cdot)>0 , \;\; \omega_Z(\cdot,J\cdot) \geq 0,
$$
therefore
$$
\Omega(\cdot,J\cdot)>0,
$$
which says that $\Omega$ is also a symplectic structure near the $\Gamma_{r_i}^s$.
Hence $(Y,\Omega)$ is  symplectic.

We now remark that the assumption on the sign of
$\lambda_i$ is inessential: suppose that
$\lambda_1>0$; then we alter $Y$ by
replacing the neighborhood of $p_1$ by $W^{sf}_{r_1}$.
Then the construction of the symplectic structure
on this $Y$ is the same.

\subsection{Proof of corollary \ref{cor_4.2.1}}\label{sect_4.3}

This follows from remark \ref{rmk_5.1.1}. If $Y$ and $Y^f$
are a pair of flops, then one of them satisfies some equation
$$
\sum_i\lambda_i\Theta_{r_i}=0
$$
and the other one satisfies
$$
-\sum_i\lambda_i\Theta_{r_i}=0.
$$
Therefore, the symplectic structures exist on them
simultaneously.

\section{Orbifold Gromov-Witten invariants of $W^s_r$ and $W^{sf}_r$}\label{sect_6}
We first introduce the cohomology group for an orbifold in the stringy sense.
Then we compute the orbifold Gromov-Witten invariants.

>From now on, $r\geq 2$ is fixed. So we drop $r$ from $W^s_r$ and $W^{sf}_r$.

\subsection{Chen-Ruan orbifold cohomology of $W^s$ and $W^{sf}$}\label{sect_6.1}
The stringy orbifold cohomology of $W^s$ is
$$
H^\ast_{CR}(W^s)=H^{\ast}(W^s)\oplus \bigoplus_{k}\cplane[\mfp]_k
\oplus\bigoplus_k\cplane[\mfq]_k.
$$

We abuse the notation here such that $[\mfp]_k$ represents the
0-cohomology of the sector $[\mfp]_k$. On the other hand, the
grading should be treated carefully: the degree of an element in
$H^\ast(W^s)$ remains the same, however the degree of $[\mfp]_k$ is
$0+\iota([\mfp]_k)$ and the same treatment applies to $[\mfq]_k$. We
call these new classes  {\em twisted classes}.

A similar definition applies to $W^{sf}$.
$$
H^\ast_{CR}(W^{sf})=H^{\ast}(W^{sf})\oplus
\bigoplus_{k}\cplane[\mfP]_k \oplus\bigoplus_k\cplane[\mfQ]_k.
$$

\subsection{Moduli spaces $\om_{0,l,k}(W^s, d[\Gamma^s], \mathbf{x}), k\geq 1$}\label{sect_6.2}

Here
$$
\mathbf{x}=(T_1,\ldots, T_k)
$$
consists of $k$ twisted sectors in $W^s$.

By the definition in \cite{CR2}, the moduli space $\om_{0,l,k}(W^s, d[\Gamma^s], \mathbf{x})$
consists of orbifold stable holomorphic maps from genus 0 curves, on which there
are $l$ smooth marked points
and  $k$ orbifold points $y_1,\ldots, y_k$, to $W^s$ such that
\begin{itemize}
\item $y_i$ are sent to $Y_i$;
\item the isotropy group at $y_i$
is $\mathbb{Z}_{|\xi^a|}$ if $y_i=[p]_a$ (or $[q]_a$), where
$|\xi^a|$ is the order of $\xi^a$;
\item the image of the map represents the homology class $d[\Gamma^s]$.
\end{itemize}
By a genus 0 curve we mean $S^2$, or a bubble tree consisting of several $S^2$'s.
The stability is the same as in the smooth case.
\begin{remark}\label{rmk_6.2.1}
There is an extra
feature for orbifold stable holomorphic maps. That is, the nodal points on a bubble tree may
also be orbifold singular points on its component: for example, say $y$ is a nodal point
that is the intersection of two spheres $S^2_+$ and $S^2_-$; then $y$ can be a singular points,
denoted by $y_+$ and $y_-$ respectively,
on both spheres. Moreover if $y_+$ is mapped to $[p]_a$, $y_-$ has to be mapped to $[p]_{r-a}$.
\end{remark}
When we write $\M_{0,l
,k}(W^s, d[\Gamma^s],\mathbf{x})$, we mean the map whose domain is $S^2$.
Usually, we call $\om$  the compactified space of $\M$ and $\M$  the top stratum of
$\om$.
\begin{lemma}\label{lemma_6.2.1}
For $k\geq 1$, the virtual dimension
$$
\dim \om_{0,0,k}(W^s, d[\Gamma^s], \mathbf{x})<0.
$$
\end{lemma}
{\bf Proof. }We recall that the virtual dimension is given by
$$
2c_1(d[\Gamma^s])+2(n-3) + k- \sum_{i=1}^k\iota(Y_i)=k-\sum_{i=1}^k\iota(Y_i)
< k-k=0.
$$
Here we use Lemma \ref{lemma_2.3.1}. q.e.d.

\v
\begin{lemma}\label{lemma_6.2.2}
$\M_{0,0,1}(W^s,d[\Gamma^s], \mathbf{x})=\emptyset$.
\end{lemma}
{\bf Proof. }This also follows from the dimension formula: the virtual dimension of
this moduli space is a {\em rational} number. q.e.d.

\subsection{Moduli spaces $\om_{0,0,0}(W^s,d[\Gamma^s])$}\label{sect_6.3}
Convention of notations: If $k=0$, it is dropped and the moduli space is denoted by
$\om_{0,l}(W^s,d[\Gamma^s])$; if $k=l=0$, then the moduli space is denoted by
$\om_{0}(W^s,d[\Gamma^s])$.

We have shown that $\om_{0,0,k}(W^s, d[\Gamma^s], \mathbf{x})$ for $k\geq 1$ has some
nice properties, following from the dimension formula.
Now we focus on $k=0$. Although its top stratum $\M_{0}(W^s,d[\Gamma^s])$ consists
of only "smooth" maps, there may be orbifold maps in lower strata. Here, by the
smoothness of a map we mean that the
domain of the map is without orbifold singularities.
The next lemma rules out this possibility.
\begin{lemma}\label{lemma_6.3.1}
$\om_{0}(W^s, d[\Gamma^s])$ only consists of smooth maps.
\end{lemma}
{\bf Proof. }If not, suppose we have a map $f\in \om_{0}(W^s, d[\Gamma^s]) $
that consists of orbifold type nodal points in the domain. By looking at the
bubble tree, we start  searching  from the leaves to look
for the first component, say $S^2_i$, that containing
a singular  nodal point. This component {\em must} contain only {\em one} singular point.
So $f|_{S^2_i}$ is an element in some moduli space $\M_{0,0,1}(W^s,d[\Gamma^s], \mathbf{x})$.
But it is claimed in Lemma \ref{lemma_6.2.2} that such an element does not exist. This proves the lemma.
q.e.d.

\v
Notice that $W^s=\tilde W^s/\mu_r$ and $\Gamma^s=\tilde\Gamma^s/\mu_r$. We may like to compare the
moduli space $\om_{0}(W^s,d[\Gamma^s])$ with $\om_{0}(\tilde W^s,d[\tilde \Gamma^s])$.
Note that $\mu_r$ acts naturally  on the latter space. We claim that
\begin{prop}\label{prop_6.3.1}
$\om_{0}(W^{s},d[\Gamma^s])=\emptyset$ if $r\nmid d$. Otherwise,
there is a natural isomorphism
$$\om_{0}(W^s,mr[\Gamma^s])=\om_{0}(\tilde W^s,m[\tilde \Gamma^s])/\mu_r.$$
if $d=mr$.
\end{prop}
{\bf Proof. }Since
$$
\om_{0}(W^s,d[\Gamma^s])=\om_0(\Gamma^s, d[\Gamma^s])
$$
and
$$
\om_{0}(\tilde W^s,d[\tilde\Gamma^s])=\om_0(\tilde \Gamma^s, d[\tilde\Gamma^s]),
$$
it is sufficient to show that
$\om_{0}(W^s, d[\Gamma^s])=\emptyset$ if $r\nmid d$ and
$$
\om_0(\Gamma^s, mr[\Gamma^s])=\om_0(\tilde \Gamma^s, m[\tilde\Gamma^s])/\mu_r.
$$
We need the following lemma.
Let $\pi: \tilde \Gamma^s\to \Gamma^s$ be the projection given by the quotient of $\mu_r$.
We claim that
\begin{lemma}\label{lemma_6.3.2}
for any
smooth map
$$
f: S^2\to \Gamma^s
$$
there is a lifting $\tilde f: S^2\to \tilde \Gamma^s$ such that $\tilde\Pi(\tilde f)=f$.
\end{lemma}
Now suppose the lemma is proved. Then we have that
$$
\om_{0}(W^{s},d[\Gamma^s])=\emptyset
$$ 
for $r\nmid d$.

To prove the second statement, we
define a map:
$$
\tilde\Pi:\om_0(\tilde \Gamma^s, m[\tilde\Gamma^s])\to
\om_0(\Gamma^s, mr[\Gamma^s])
$$
given by $\tilde\Pi(\tilde f)=\pi\circ\tilde f$. It is clear that this induces an injective map
$$
\Pi:\om_0(\tilde \Gamma^s, m[\tilde\Gamma^s])/\mu_r\to
\om_0(\Gamma^s, mr[\Gamma^s]).
$$
On the other hand,  since
a stable smooth map on a bubble tree consists of  smooth maps on each component
of the tree that match at nodal points,
therefore, by Lemma \ref{lemma_6.3.2} the map can be components wise lifted. This shows that
$\Pi$ is surjective. q.e.d.
\v\n
{\bf Proof of  Lemma \ref{lemma_6.3.2}: }
$S^2$ and $\Gamma^s$ are $\pone^1$. We identify them as $\cplane\cup \{\infty\}$ as usual.
On $\Gamma^s$, we assume $\mfp$ and $\mfq$ are $0$ and $\infty$ respectively.

Suppose that
$$
\Lambda_0=f^{-1}(\mfp)=\{ x_1,\ldots,x_m\}, \;\;,
\Lambda_\infty=f^{-1}(\mfq)=\{y_1,\ldots,y_n\}.
$$
 Let $z$ be the
coordinate of the first sphere; we write
$$
f(z)=[p(z),q(z)].
$$
Now since $f$ is assumed to be smooth at the $x_i$, the map can be lifted
with respect to the uniformizing system of $\mfp$:  namely, suppose that
$$
\pi_\mfp: D_{\epsilon}(0)\subset\cplane \to D_{\epsilon^r}(\mfp)\cplane; \;\; \pi_\mfp(w)=w^r
$$
gives the uniformizing system of the neighborhood of $\mfp$ for some
$\epsilon$;  $f$, restricted to a small neighborhood $U_{x_i}$,
 can be lifted
to
$$
\tilde f: U_{x_i}\to D_\epsilon
$$
such that $f=\pi_\mfp\circ\tilde f$. Without loss of generality, we assume that $f(U_{x_i})
=D_\epsilon(0)$.
 Therefore we have a lifting
$$
\tilde f: \bigcup_{i}U_{x_i}\cup \bigcup_j U_{y_j}\to D_{\epsilon}(0)\cup D_{\epsilon}(\infty)
$$
for $f$. Now we look at the rest of the map
$$
f: S^2-\bigcup_{i}U_{x_i}\cup \bigcup_j U_{y_j}\to \Gamma^s-D_{\epsilon^r}(\mfp)
\cup D_{\epsilon^r}(\mfq).
$$
We ask if this map can be lifted to the covering space
$$
\tilde \Gamma^s-D_{\epsilon^r}(0)
\cup D_{\epsilon^r}(\infty)\to  \Gamma^s-D_{\epsilon^r}(\mfp)
\cup D_{\epsilon^r}(\mfq).
$$
The answer is affirmative by the elementary lifting theory for the covering space. Therefore,
the whole map $f$ has a lifting $\tilde f$. The ambiguity of the lifting is up to the
$\mu_r$ action.  q.e.d.

\subsection{Orbifold Gromov-Witten invariants on $W^s$}\label{sect_6.4}

We study the Gromov-Witten invariants that are needed in this paper.

Given a moduli space $\om_{0,l,k}(W^s,d[\Gamma^s], \mathbf{x})$, one can define
the Gromov-Witten invariants via evaluation maps:
\begin{eqnarray*}
&&ev_i:\om_{0,l,k}(W^s,d[\Gamma^s], \mathbf{x})\to X, 1\leq i\leq l;\\
&&{ev}^{orb}_j:\om_{0,l,k}(W^s,d[\Gamma^s], \mathbf{x})\to Y_j, 1\leq j\leq k.
\end{eqnarray*}
The Gromov-Witten invariants are given by
\begin{eqnarray*}
&&\Psi_{(d[\Gamma^s],0,l,k,\mathbf{x})}^{W^s}(\alpha_1,\ldots, \alpha_l,\gamma_1,\ldots,\gamma_k)
\\
&&=
\int_{[\om_{0,l,k}(W^s,d[\Gamma^s], \mathbf{x})]^{vir}}
\bigcup_i ev_i^\ast(\alpha_i)\cup\bigcup_j {ev}_j^{orb,\ast}(\beta_j).
\end{eqnarray*}
Here $\alpha_i\in H^\ast(X)$ and $\beta_j\in H^\ast(Y_j)$.
Note that $l,k$ and $\mathbf{x}$ are specified by the $\alpha_i$ and $\beta_j$.
For the sake of simplicity and consistency, we  also re-denote the invariants by
$$
\Psi_{(d[\Gamma^s],0,l+k)}^{W^s}(\alpha_1,\ldots, \alpha_l,\gamma_1,\ldots,\gamma_k),
$$
when the $\alpha_i$ and $\beta_j$ are given.
 \begin{lemma}\label{lemma_6.4.1}
For $k\geq 1$ and $d\geq 1$
$$
\Psi_{(d[\Gamma^s],0,0,k,\mathbf{x})}^{W^s}=0.
$$
\end{lemma}
{\bf Proof. }As explained in Lemma \ref{lemma_6.2.1}, this moduli space has negative
dimension. Therefore the Gromov-Witten invariants have to be 0. q.e.d.
\v\n
\begin{prop}\label{prop_6.4.1}
For $d\geq 1$, if $r\nmid d$,
$\Psi_{(d[\Gamma^s],0)}^{W^s}$ vanishes. Otherwise, if $d=mr$
$$
\Psi_{(mr[\Gamma^s],0)}^{W^s}=\frac{1}{m^3}.
$$
\end{prop}
{\bf Proof. }We have shown that
$$
\om_{0}(W^s,mr[\Gamma^s])=\om_{0}(\tilde W^s,m[\tilde\Gamma^s])/\mu_r.
$$
This would suggest that
\begin{equation}\label{eqn_6.1}
\Psi_{(mr[\Gamma^s],0)}^{W^s}=\frac{1}{r} \Psi_{(m[\tilde \Gamma^s],0)}^{\tilde W^s}.
\end{equation}
This has to be shown by virtual techniques. Following the standard  construction
of virtual neighborhoods of moduli spaces, we have a smooth virtual moduli space
$$
\mc U_{0}(\tilde W^s,m[\tilde\Gamma^s])\supset \om_{0}(\tilde W^s, m[\tilde\Gamma^s]),
$$
with an obstruction bundle $\tilde{\mc O}$. The Gromov-Witten invariant is then given by
$$
\Psi_{(m[\tilde \Gamma^s],0)}^{\tilde W^s}
=\int_{\mc U_{0}(\tilde W^s,m[\tilde\Gamma^s])}\Theta(\tilde{\mc O}).
$$
Here
$\Theta(\tilde {\mc O})$ is the Thom form of the bundle.
See the construction of virtual neighborhood
in \cite{CL} (and orginally
in \cite{R2}).
The construction of virtual neighborhoods for $\om_{0}(W^s, mr[\Gamma^s])$ is parallel.
We also have
$$\mc U_{0}( W^s,mr[\Gamma^s])$$
with obstruction bundle $\mc O$.
The model can be suggestively expressed as
$$
(\mc U_{0}(W^s,mr[\Gamma^s]), {\mc O})
=(\mc U_{0}(\tilde W^s,m[\tilde\Gamma^s]), \tilde{\mc O})/\mu_r.
$$
Therefore, we conclude that
$$
\Psi_{(mr[\Gamma^s],0)}^{W^s}=\frac{1}{r}\int_{\mc U_{0}
(\tilde W^s,m[\tilde\Gamma^s])}\Theta(\tilde{\mc O})=\frac{1}{r}\Psi_{(m[\tilde \Gamma^s],0)}^{\tilde W^s}.
$$

On the other hand,
$$
\Psi_{(m[\tilde \Gamma^s],0)}^{\tilde W^s}=\frac{r}{m^3}.
$$
This is computed in \cite{BLK}. Therefore the proposition is proved. q.e.d.

\subsection{3-point functions on  $H^\ast_{CR}(W^s)$ and
$H^\ast_{CR}(W^{sf})$}\label{sect_6.5}

On $W^s$,
$$
H^\ast_{CR}(W^s)=\cplane[1]\oplus \cplane(H^s)\oplus \bigoplus_{i=1}^{r-1}\cplane[\mfp]_i
\oplus\bigoplus_{j=1}^{r-1}\cplane[\mfq]_j.
$$

Given $\beta_i,1\leq i\leq 3,$ in $H^\ast_{CR}(X)$ one  defines  the 3-point
function
as following:
$$
\Psi^{W^s}(\beta_1,\beta_2,\beta_3)
=\Psi^{W^s}_{CR}(\beta_1,\beta_2,\beta_3)+\sum_{d\geq 1}
\Psi_{(d[\Gamma^s],0,3)}^{W^s}(\beta_1,\beta_2,\beta_3)q^{d[\Gamma^s]}.
$$
 Here the first term
$$
\Psi^{W^s}_{CR}(\beta_1,\beta_2,\beta_3)=\Psi_{([0],0,3)}^{W^s}(\beta_1,\beta_2,\beta_3)
$$
is the 3-point function
defining the Chen-Ruan product. In the smooth case, this is just
$$
\int\beta_1\wedge\beta_2\wedge\beta_3.
$$
A similar expression for the orbifold case still holds. This is
proved  in \cite{CH}: by introducing twisting factors, one can turn a twisted
form $\beta$ on twisted sector into a formal form $\tilde \beta$ on the global orbifold.
Then we still have
$$
\Psi^{W^s}_{CR}(\beta_1,\beta_2,\beta_3)
=\int^{orb}_{W^s}\tilde\beta_1\wedge\tilde\beta_2\wedge\tilde\beta_3.
$$
\begin{remark}\label{rmk_6.5.1}
Unfortunately, for the local model, $\Psi^{W^s}_{cr}(\beta_1,\beta_2,
\beta_3)$ does not make sense if and only if  all  $\beta_i$ are smooth classes,
for the moduli space of the latter case is non-compact. Hence $\Psi^{W^s}_{CR}(\beta_1,\beta_2,
\beta_3)$ is only a notation at the moment. But we will need it when we move on to study
compact symplectic conifolds.
\end{remark}

By the computation in \S \ref{sect_6.4},
we have
\begin{prop}\label{prop_6.5.1}
If at least one of the $\beta_i$ is a twisted class,
$$
\Psi^{W^s}(\beta_1,\beta_2,\beta_3)
=\Psi^{W^s}_{CR}(\beta_1,\beta_2,\beta_3).
$$
\end{prop}
{\bf Proof. }Case 1, if all $\beta_i$ are twisted classes,
$$\Psi_{( d[\Gamma^s],0,3)}^{W^s}(\beta_1,\beta_2,\beta_3)=0$$
if $d\geq 1$.
\v
Now suppose $\beta_3$ is not twisted and the other two are.
\v\n
Case 2: Suppose $\beta_3=1$; then it is well known that
$$\Psi_{( d[\Gamma^s],0,3)}^{W^s}(\beta_1,\beta_2,1)=0$$
if $d\geq 1$.

\v\n
Case 3: suppose that   $\beta_3=nH^s$; then
$$
\Psi_{( d[\Gamma^s],0,3)}^{W^s}(\beta_1,\beta_2,\beta_3)
= \beta_3(d[\Gamma^s])
\Psi_{(d[\Gamma^s],0,2)}^{W^s}(\beta_1,\beta_2)=0.
$$
\v
Similar arguments can be applied to the case in which  only one of the $\beta_i$ is twisted.
Hence the claim follows. q.e.d.

\v
Now suppose $\deg(\beta_i)=2$, i.e. $\beta_i=n_i H^s$. Then

$$
\sum_{m\geq 1}
\Psi_{ (mr[\Gamma^s],0,3)}^{W^s}(\beta_1,\beta_2,\beta_3)q^{mr[\Gamma_s]}
=\beta_1([r\Gamma^s])\beta_2([r\Gamma^s])\beta_3([r\Gamma^s])\frac{q^{[r\Gamma^s]}}{1-q^{[r\Gamma^s]}}.
$$
The last statement follows from Proposition \ref{prop_6.4.1}.
Hence
$$
\Psi^{W^s}(\beta_1,\beta_2,\beta_3)
=\int_{W^s}^{orb}\beta_1\wedge\beta_2\wedge\beta_3+\beta_1([r\Gamma^s])\beta_2([r\Gamma^s])\beta_3([r\Gamma^s])
\frac{q^{[r\Gamma^s]}}{1-q^{[r\Gamma^s]}}.
$$
Formally, we write $[\tilde\Gamma^s]=[r\Gamma_s]$.
To summarize,
\begin{prop}\label{prop_6.5.2}
The three-point function
$
\Psi^{W^s}(\beta_1,\beta_2,\beta_3) $ of $W^s$
is
$$
\Psi^{W^s}_{CR}(\beta_1,\beta_2,\beta_3)
$$
if
 at least one of the $\beta_i$ is twisted or of degree 0, or
$$
\Psi^{W^s}_{cr}(\beta_1,\beta_2,\beta_3)
+\beta_1(\tilde\Gamma^s)\beta_2(\tilde\Gamma^s)\beta_3(\tilde\Gamma^s)
\frac{q^{[\tilde\Gamma^s]}}{1-q^{[\tilde\Gamma^s]}},
$$
if $\deg(\beta_i)=2, 1\leq i\leq 3.$
\end{prop}
This proposition says that
the quantum product $\beta_1\star\beta_2$ is the usual product( in the sense of the Chen-Ruan
ring structure) except for the case in which $\deg(\beta_1)=\deg(\beta_2)=2$.
Next, we write down the Chen-Ruan ring structure for twisted classes:

\begin{prop}\label{prop_6.5.3}
The Chen-Ruan products for twisted classes are given by
\begin{eqnarray*}
&& [\mfp]_i\star[\mfq]_j=0, \\
&& [\mfp]_i\star[\mfp]_j=\delta_{i+j,r}\Theta_{\mathbf p},\\
&& [\mfq]_i\star[\mfq]_j=\delta_{i+j,r}\Theta_{\mathbf q}.
\end{eqnarray*}
Here $\Theta_p$ and $\Theta_q$ are Thom forms of the normal bundles of $\mathbf{p}$ and
$\mathbf q$
in $W^s$.
Also
$$
\beta\star H^s=0
$$
if $\beta$ is a twisted class.
\end{prop}
{\bf Proof. }This follows from the theorem in \cite{CH}. As an example,  we verify
$$
[\mfp]_i\star[\mfp]_j=\delta_{i+j,r}\Theta_{\mathbf p}=0.
$$
For other cases, the proof is similar.
The normal bundle of $\mathbf p$ is a rank 3  orbi-bundle which splits as three lines
$\cplane_p, \cplane_y$ and $\cplane_z$ (cf. S\ref{sect_1.25}). Let $\Theta_p, \Theta_y$ and
$\Theta_z$ be the corresponding Thom forms. Then the twisting factor(cf. \cite{CH})
 of $[\mfp]_i$
is
$$
\mk t([\mfp]_i)=
\Theta_p^{b}\Theta_y^{r-b}\Theta_z^i.
$$
Here $b\equiv ai (\mod r)$
is an integer
between $0$ and $r-1$.
Similarly, we write
$$
\mk t([\mfp]_j)=
\Theta_p^{c}\Theta_y^{r-c}\Theta_z^j.
$$
Here $c\equiv aj(\mod r)$
is an integer
between $0$ and $r-1$.
Then we have a formal computation
$$
[\mfp]_i\star[\mfp]_j
=\mk t([\mfp]_i)\wedge \mk t([\mfp]_j)
=\delta_{i+j,r}\Theta_{\mathbf{p}}.
$$
q.e.d.

\v
Equivalently, this can be restated in terms of $\Psi^{W^s}_{cr}$ as
\begin{prop}\label{prop_6.5.4}
Suppose at least one of the
$\beta_i$ is twisted in the three-point function $\Psi^{W^s}_{cr}(\beta_1,\beta_2,\beta_3)$.
Then only the following functions are nontrivial:
\begin{eqnarray*}
\Psi^{W^s}_{cr}([\mfp]_i,[\mfp]_j,1)&=&\delta_{i+j,r}\frac{1}{r};\\
\Psi^{W^s}_{cr}([\mfq]_i,[\mfq]_j,1)&=&\delta_{i+j,r}\frac{1}{r}.
\end{eqnarray*}
\end{prop}

\subsection{Identification of three-point functions $\Psi^{W^s}$ and $\Psi^{W^{sf}}$}\label{sect_6.6}

We follow the argument in \cite{LR}. Define a map
$$
\phi: H^\ast_{CR}(W^{sf})\to H^\ast_{CR}(W^{s}).
$$
On twisted classes, we define
$$
\phi([\mfP]_k)=[\mfp]_k,
\;\;
\phi([\mfQ]_k)=[\mfq]_k.
$$
And on $H^\ast_{CR}(W^{sf})$, $\phi$ is defined as in the smooth
case in \cite{LR}. Since at the moment we are working in the local model,
we should avoid using  Poincare duality. We give a direct
geometric construction of the map. On the other hand, a technical
issue mentioned in Remark \ref{rmk_6.5.1} is dealt with: let $\beta_i^{sf}, 1\leq
i\leq 3,$ be 2-forms on $W^{sf}$ representing the classes
$[\beta_i^{sf}]$; by the identification of $W^{sf}-\Gamma^{sf}$
with $W^s-\Gamma^s$, we then also have 2-forms in $W^s-\Gamma^s$
which as cohomology classes can be {\em uniquely} extended over
$W^s$. The cohomology classes are denoted by
$$
[\alpha_i]=\phi([\beta_i]).
$$
Moreover we can require that the representing forms, denoted by
$
\alpha_i,$
coincide with $\beta_i$ away from  the $\Gamma$'s.

Then we can define
\begin{eqnarray*}
\Psi^{W^s}_{CR}([\alpha_1],[\alpha_2],[\alpha_3])
&-&\Psi^{W^{sf}}_{CR}([\beta_1],[\beta_2],[\beta_3])\\
&:=&\int_{W^s}^{orb}\alpha_1\wedge\alpha_2\wedge\alpha_3
-\int_{W^{sf}}^{orb}\beta_1\wedge\beta_2\wedge\beta_3.
\end{eqnarray*}
The well-definedness can be easily seen due to the coincidence of the $\alpha_i$ and
$\beta_i$  outside a compact set.
Moreover,
\begin{lemma}\label{6.6.1}
Suppose that $\deg \beta_i=2$; then
\begin{eqnarray*}
\Psi^{W^s}_{CR}([\alpha_1],[\alpha_2],[\alpha_3])
-\Psi^{W^{sf}}_{CR}([\beta_1],[\beta_2],[\beta_3])
&=&\alpha_1(\tilde\Gamma^s)\alpha_2(\tilde\Gamma^s)\alpha_3(\tilde\Gamma^s)\\
&=&-\beta_1(\tilde\Gamma^{sf})\beta_2(\tilde\Gamma^{sf})\beta_3(\tilde\Gamma^{sf}).
\end{eqnarray*}
\end{lemma}
{\bf Proof. }We lift the problem to  $\tilde W^s$ and $\tilde W^{sf}$. Then
we can further  deform both models simultaneously to $\tilde V^s$
and $\tilde V^{sf}$ as \cite{F}.
Each of them consists $r$ copies of the standard model $\mc O(-1)\oplus \mc O(-1)
\to \pone^1$.
$\tilde V^{sf}$
is a flop of $\tilde V^s$. Therefore,  the computations are essentially
$r$ copies of
the computation on the standard model. By the argument in \cite{LR}, we have
\begin{eqnarray*}
\int_{W^s}^{orb}\alpha_1\wedge\alpha_2\wedge\alpha_3
&-&\int_{W^{sf}}^{orb}\beta_1\wedge\beta_2\wedge\beta_3\\
&=&
\frac{1}{r}\left(
\int_{\tilde W^s}\alpha_1\wedge\alpha_2\wedge\alpha_3
-\int_{\tilde W^{sf}}\beta_1\wedge\beta_2\wedge\beta_3\right)\\
&=&
\frac{1}{r}\cdot r\cdot\alpha_1(\tilde\Gamma^s)\alpha_2(\tilde\Gamma^s)\alpha_3
(\tilde\Gamma^s)\\
&=&\alpha_1(\tilde\Gamma^s)\alpha_2(\tilde\Gamma^s)\alpha_3
(\tilde\Gamma^s).
\end{eqnarray*}

Now we conclude that
\begin{theorem}\label{theorem_6.6.1}
Let $\beta_i\in H^\ast_{CR}(W^{sf}),1\leq i\leq 3,$ and
$\alpha_i=\phi(\beta_i)$. Then
$$
\Psi^{W^s}(\alpha_1,\alpha_2,\alpha_3)=\Psi^{W^{sf}}(\beta_1,\beta_2,\beta_3)
$$
with the identification of $[\Gamma^s]\leftrightarrow -[\Gamma^{sf}]$.
\end{theorem}
{\bf Proof. }The only nontrivial verification is for all
$\deg \beta_i=2$. Suppose this is the case. Then the difference
$$
\Psi^{W^s}(\alpha_1,\alpha_2,\alpha_3)-\Psi^{W^{sf}}(\beta_1,\beta_2,\beta_3)
$$
includes two parts. Part(I) is
$$
\Psi^{W^s}_{cr}([\alpha_1],[\alpha_2],[\alpha_3])
-\Psi^{W^{sf}}_{cr}([\beta_1],[\beta_2],[\beta_3])
=\alpha_1(\tilde\Gamma^s)\alpha_2(\tilde\Gamma^s)\alpha_3(\tilde\Gamma^s)
$$
and part(II) is
\begin{eqnarray*}
&&\alpha_1(\tilde\Gamma^s)\alpha_2(\tilde\Gamma^s)\alpha_3(\tilde\Gamma^s)
\frac{q^{[\tilde\Gamma^s]}}{1-q^{[\tilde\Gamma^s]}}
-\beta_1(\tilde\Gamma^{sf})\beta_2(\tilde\Gamma^{sf})\beta_3(\tilde\Gamma^{sf})
\frac{q^{[\tilde\Gamma^{sf}]}}{1-q^{[\tilde\Gamma^{sf}]}}\\
&=&\alpha_1(\tilde\Gamma^s)\alpha_2(\tilde\Gamma^s)\alpha_3(\tilde\Gamma^s)
\frac{q^{[\tilde\Gamma^s]}}{1-q^{[\tilde\Gamma^s]}}
+\alpha_1(\tilde\Gamma^s)\alpha_2(\tilde\Gamma^s)\alpha_3(\tilde\Gamma^s)
\frac{q^{[-\tilde\Gamma^{s}]}}{1-q^{[-\tilde\Gamma^{s}]}}\\
&=&-\alpha_1(\tilde\Gamma^s)\alpha_2(\tilde\Gamma^s)\alpha_3(\tilde\Gamma^s).
\end{eqnarray*}
Here we use $[\Gamma^s]\leftrightarrow -[\Gamma^{sf}]$.
Part(I) cancels part (II), therefore
$$
\Psi^{W^s}(\alpha_1,\alpha_2,\alpha_3)=\Psi^{W^{sf}}(\beta_1,\beta_2,\beta_3).
$$
q.e.d.
\v

\section{Ruan's conjecture on orbifold symplectic flops}\label{sect_7}

\subsection{Ruan cohomology}\label{sect_7.1}
Let $X$ and $Y$ be compact symplectic orbifolds related by
symplectic flops. Correspondingly, $\Gamma_i^s$ and
$\Gamma_i^{sf}$, $1\leq i\leq k$, are extremal rays on $X$ and $Y$
respectively.  We define
 three-point functions on $X$ (similarly on $Y$):
$$
\Psi^X_{qc}(\beta_1,\beta_2,\beta_3)
=\Psi^X_{CR}(\beta_1,\beta_2,\beta_3)
+\sum_{i=1}^k \sum_{d=1}^\infty
\Psi^X_{(d[\Gamma^s_i], 0, 3)}(\beta_1,\beta_2,\beta_3).
$$
This induces a ring structure on $H^\ast_{CR}(X)$
\begin{defn}\label{defn_7.1}
Define the product on $H^\ast_{CR}(X)$ by
$$
\langle\beta_1\star_r\beta_2,\beta_3\rangle
=\Psi^X_{qc}(\beta_1,\beta_2,\beta_3).
$$
We call this the Ruan product on $X$. This cohomology ring is
denoted by $RH_{CR}(X)$.
\end{defn}
Similarly,  we can define $RH^\ast_{CR}(Y)$ by using the three-point
functions given by $\Gamma_i^{sf}$. Ruan conjectures that
\begin{conjecture}[Ruan]
$RH^\ast_{CR}(X)$ is isomorphic to $RH^\ast_{CR}(Y)$.
\end{conjecture}

\subsection{Verification of Ruan's conjecture}\label{sect_7.2}
Set
$$
\Phi([\Gamma^s_u])=-[\Gamma^{sf}_u].
$$
This induces an obvious identification
$$
\Phi:H_2(X)\to H_2(Y).
$$
As explained in the local model, there is a natural isomorphism
$$
\phi: H^\ast_{CR}(Y)\to H^\ast_{CR}(X).
$$
We explain $\phi$.
For twisted classes $[\mfP_s]_i$ and $[\mfQ_t]_j$ we define
$$
\phi([\mfP_u]_i)=[\mfp_u]_i, \;\;,
\phi([\mfQ_v]_j)=[\mfq_v]_j.
$$
For degree $0$ or $6$-forms, $\phi$ is defined in an obvious way.
For $\alpha\in H^{2}_{orb}(Y)$, $\phi(\alpha)$ is defined to be the unique
extension of
$$
\alpha|_{X-\cup\Gamma_u^s}=\alpha|_{Y-\cup\Gamma_v^{sf}}
$$
over $X$. For  $\beta\in H^4(Y)$, define $\phi(\beta)\in H^4(X)$ to be
the extension  as above such that
$$
\int_{X}\phi(\beta)\wedge\phi(\alpha)=\int_Y\beta\wedge\alpha,
$$
for any $\alpha\in H^2(Y)$.
Then
\begin{theorem}\label{theorem_7.2.1}
For any classes $\beta_i\in H^\ast_{CR}(Y), 1\leq i\leq 3,$
$$
\Phi_\ast(\Psi^X_{qc,r}(\phi(\beta_1),\phi(\beta_2),\phi(\beta_3)))
=\Psi^Y_{qc,r}(\alpha_1,\alpha_2,\alpha_3).
$$
\end{theorem}
{\bf Proof. }If one of $\beta_i$, say $\beta_1$, has degree $\geq 4$, the quantum correction
term vanishes. Therefore, we need only verify
$$
\Psi^X_{CR}(\phi(\beta_1),\phi(\beta_2),\phi(\beta_3))
=\Psi^Y_{CR}(\alpha_1,\alpha_2,\alpha_3).
$$
We  choose $\beta_1$ to be supported away from the $\Gamma^{sf}$. Then we have following observations:
\begin{itemize}
\item whenever $\beta_2$ or $\beta_3$ is a twisted class, both sides are equal to 0;
\item if $\beta_2$ and $\beta_3$ are in $H^\ast(Y)$, then
\begin{eqnarray*}
\Psi^X_{cr}(\phi(\beta_1),\phi(\beta_2),\phi(\beta_3))
&=&\int_X\phi(\beta_1)\wedge\phi(\beta_2)\wedge\phi(\beta_3)\\
&=&\int_Y\beta_1\wedge\beta_2\wedge\beta_3
=\Psi^Y_{cr}(\alpha_1,\alpha_2,\alpha_3).
\end{eqnarray*}
\end{itemize}

Now we assume that $\beta_i$ are either twisted classes or degree 2 classes. Then the verification
is exactly same as that in Theorem \ref{theorem_6.6.1}. q.e.d.

\v\n
As an corollary, we have proved
\begin{theorem}\label{theorem_7.2.2}
Suppose $X$ and $Y$ are related via an orbifold symplectic flops,
Via the map $\phi$ and coordinate change $\Phi$,
$$
RH^\ast_{CR}(X)\cong RH^\ast_{CR}(Y).
$$
\end{theorem}
This explicitly realizes the claim of Theorem \ref{theorem_1.2}.


\begin{thebibliography}{L3}

\bibitem[BKL]{BLK} J. Bryant,  S. Katz,  N. Leung, Multiple covers and the integrality conjecture for rational curves in
CY threefolds , J. ALgebraic Geometry 10(2001),no.3.,549-568.
\bibitem[CH]{CH} B. Chen, S. Hu, A de Rham model of Chen-Ruan cohomology ring of
 abelian orbifolds, To appear in Math. Ann.
\bibitem[CL]{CL} B. Chen, A-M. Li, Symplectic Virtual Localization of Gromov-Witten invariants,
                in preparation.
\bibitem[CLZZ]{CLZZ}  B. Chen, A-M. Li, Q. Zhang, G. Zhao, Relative Gromov-Witten invariants
                  and glue formula, in preparation.
\bibitem{CT}{CT} B. Chen, G. Tian, Virtual orbifolds and Localization, preprint.
\bibitem[CR1]{CR1}  W. Chen, Y. Ruan, A new cohomology theory for
                  orbifold, AG/0004129, Commun. Math. Phys.,
                  248(2004), 1-31.
\bibitem[CR2]{CR2}  W. Chen, Y. Ruan, orbifold Gromov-Witten
                   theory,  AG/0011149. Cont. Math., 310, 25-86.
\bibitem[CR3]{CR3}  W. Chen, Y. Ruan, orbifold quantum cohomology,  Preprint AG/0005198.

\bibitem[F]{F} R. Friedman, Simultaneous resolutions of threefold double points, Math. Ann.
274(1986) 671-689.
\bibitem[Gr]{Gr} M. Gromov, Pseudo holomorphic curves in symplectic manifolds,
                Invent. math., 82 (1985), 307-347.
\bibitem[HZ]{HZ} J. Hu, W. Zhang, Mukai flop and Ruan cohomology,
             Math. Ann. 330, No.3, 577-599 (2004).

\bibitem[K]{K} J. Koll\'ar, Flips, Flops, Minimal Models, Etc.,
                  Surveys in Differential Geometry, 1(1991),113-199.

\bibitem[La]{Laufer}  Henry B. Laufer, On $CP^1$ as an xceptional
            set, In recent developments in several complex variables
            ,261-275, Ann. of Math. Studies 100, Princeton, 1981.

\bibitem[L]{L} E. Lerman, Symplectic cuts, Math Research Let 2(1985) 247-258

 \bibitem[LR]{LR}  A-M. Li, Y. Ruan, Symplectic surgery and
                     Gromov-Witten invariants of Calabi-Yau 3-folds, Invent. Math. 145,
                   151-218(2001)
\bibitem[LZZ]{LZZ} A-M. Li, G. Zhao, Q. Zheng, The number of ramified covering of a Riemann surface by Riemann
surface, Commu. Math. Phys, 213(2000), 3, 685--696.
\bibitem[Reid]{Reid}  M. Reid, Young Person's Guide to Canonical
                      Singularities, Proceedings of Symposia in Pure Mathematics, V.46 (1987).
\bibitem[R]{R1} Y. Ruan, Surgery, quantum cohomology and birational geometry, math.AG/9810039.

\bibitem[R2]{R2} Y. Ruan, Virtual neighborhoods and pseudo-holomorphic curves, alg-geom/9611021 ,


\bibitem[S]{S}  I. Satake, The Gauss-Bonnet theorem for
                  V-manifolds, J. Math. Soc. Japan 9(1957), 464-492.

\bibitem[STY]{STY}  I. Smith, R.P. Thomas, S.-T. Yau, Symplectic
                conifold transitions.   SG/0209319. J. Diff.
                Geom., 62(2002), 209-232.

\end{thebibliography}
\end{document}